\documentclass[preprint,times,12pt]{elsarticle}

\usepackage{graphicx}
\usepackage{dcolumn}
\usepackage{bm}
\usepackage{xcolor}
\usepackage{subcaption}
\usepackage{amsmath,amssymb,amsfonts}
\usepackage{hyperref}
\hypersetup{colorlinks=true,
            linkcolor=blue,
            anchorcolor=blue,
            citecolor=blue}
\usepackage{natbib} \setcitestyle{comma,sort&compress,numbers,square}
\usepackage{amsthm}
\usepackage{geometry}
\usepackage{algpseudocode}
\usepackage{algorithm}
\theoremstyle{definition}

\usepackage{cleveref}
\usepackage{comment}
\usepackage{multirow}
\usepackage{xfrac}
\usepackage{booktabs}
\usepackage[utf8]{inputenc}

\makeatletter
\def\ps@pprintTitle{%
 \let\@oddhead\@empty
 \let\@evenhead\@empty
 \def\@oddfoot{}%
 \let\@evenfoot\@oddfoot}
\makeatother

\allowdisplaybreaks

\newcommand{\Break}{\State \textbf{break} }

\DeclareMathOperator*{\argmin}{arg\,min}

\usepackage{xcolor}

\definecolor{mat5}{rgb}{0.0, 0.0, 0.5156}
\definecolor{mat4}{rgb}{0.0, 0.3125, 1.0}
\definecolor{mat3}{rgb}{0.1094, 1.0000, 0.8906}
\definecolor{mat2}{rgb}{0.9062, 1.0000, 0.0938}
\definecolor{mat1}{rgb}{1.0, 0.2969, 0}

\definecolor{mc_mat1}{rgb}{0.5, 1.0, 0.5}
\definecolor{mc_mat2}{rgb}{0, 0, 0.5156}

\newcounter{rowcntr}[table]
\renewcommand{\therowcntr}{(\alph{rowcntr})}

\newcommand{\customref}[2]{\ref{#1}\ref{#2}}

\newcolumntype{N}{>{\refstepcounter{rowcntr}\therowcntr}c}

\AtBeginEnvironment{tabular}{\setcounter{rowcntr}{0}}

\usepackage{xcolor}
\usepackage{pict2e}

\newsavebox{\ORCIDlogo}
\savebox{\ORCIDlogo}{%
\setlength{\unitlength}{\dimexpr 1em/256\relax}%
\begin{picture}(256,256)%
  \color[HTML]{A6CE39}\put(128,128){\circle*{256}}%
  \color{white}%
  \put(78.6,199.2){\circle*{20}}%
  \moveto(70.9,176,9)\lineto(86.3,176,9)\lineto(86.3,69.8)\lineto(70.9,69.8)%
  \closepath\fillpath%
  \moveto(108.9,176.9)\lineto(150.5,176.9)%
  \curveto(190.1,176.9)(207.5,148.6)(207.5 ,123.3)%
  \curveto(207.5,95,8)(186,69.7)(150.7,69.7)%
  \lineto(108.9,69.7)%
  \closepath\fillpath%
  \color[HTML]{A6CE39}%
  \moveto(124.3,83.6)\lineto(148.8,83.6)%
  \curveto(183.7,83.6)(191.7,110.1)(191.7,123.3)%
  \curveto(191.7,144.8)(178,163)(148,163)%
  \lineto(124.3,163)%
  \closepath\fillpath%
\end{picture}%
}

\newcommand\orcidicon[1]{\href{https://orcid.org/#1}{\usebox{\ORCIDlogo}}}

\begin{document}
\begin{frontmatter}
\title{Discrete Variable Topology Optimization Using Multi-Cut Formulation and Adaptive Trust Regions}
\date{}
\author[WISC]{Zisheng Ye \orcidicon{0000-0001-6675-9747}}
\author[WISC]{Wenxiao Pan \orcidicon{0000-0002-2252-7888} \corref{cor}}
\ead{wpan9@wisc.edu}
\cortext[cor]{Corresponding author}
\address[WISC]{Department of Mechanical Engineering, University of Wisconsin-Madison, Madison, WI 53706, USA}

\begin{abstract}
    We present a new framework for efficiently solving general topology optimization (TO) problems that find an optimal material distribution within a design space to maximize the performance of a part or structure while satisfying design constraints. These problems can involve convex or non-convex objective functions and may include multiple candidate materials. The framework is designed to greatly enhance computational efficiency, primarily by diminishing optimization iteration counts and thereby reducing the frequency of solving associated state-equilibrium partial differential equations (PDEs) (e.g., through the finite element method (FEM)). It maintains binary design variables and addresses the large-scale mixed integer nonlinear programming (MINLP) problem that arises from discretizing the design space and PDEs. The core of this framework is the integration of the generalized Benders' decomposition and adaptive trust regions. Specifically, by formulating the master sub-problem (decomposed from the MINLP) as a multi-cut optimization problem and enabling the estimation of the upper and lower bounds of the original objective function, significant acceleration in solution convergence is achieved. The trust-region radius adapts based on a merit function. To mitigate ill-conditioning due to extreme parameter values, we further introduce a parameter relaxation scheme where two parameters are relaxed in stages at different paces, improving both solution quality and efficiency. Numerical tests validate the framework's superior performance, including minimum compliance and compliant mechanism problems in single-material and multi-material designs. We compare our results with those of other methods and demonstrate significant reductions in optimization iterations (and thereby the number of FEM analyses required) by about one order of magnitude, while maintaining comparable optimal objective function values and material layouts. As the design variables and constraints increase, the framework maintains consistent solution quality and efficiency, underscoring its good scalability. We anticipate this framework will be especially advantageous for TO applications involving substantial design variables and constraints and requiring significant computational resources for FEM analyses (or PDE solving).
\end{abstract}


\end{frontmatter}

\section{Introduction}

Topology optimization (TO) seeks an optimal material distribution within a given design space to maximize the performance of a part or structure, while satisfying specific design constraints such as volume or material usage. TO has seen growth in applications across a wide spectrum of fields, including multi-scale mechanics \cite{wu2021topology}, fluid dynamics \cite{borrvall2003topology, gersborg2005topology, kreissl2011topology}, electromagnetics \cite{kiziltas2003topology}, photonics \cite{nomura2007structural}, and coupled multi-physics problems \cite{dede2009multiphysics, kreissl2010topology, kumar2020topology}. It has also been widely adopted in the automotive \cite{cavazzuti2011high} and aerospace industries \cite{zhu2016topology}. The physics governing the system—whether structural mechanics, fluid dynamics, or multiphysics—can be described by state equilibrium equations, typically represented as partial differential equations (PDEs). Thus, TO can be mathematically formulated as a PDE-constrained optimization problem. The growth of TO applications has introduced significant computational challenges due to the vast and complex search space involved in optimization, the substantial computational expense of solving PDEs, and the necessity of repeatedly solving the PDEs throughout the optimization process. To address these challenges, investigating new methodologies and algorithms is imperative.

The present work concerns general and practically relevant TO problems: the objective functions can be either convex or non-convex with respect to the design variables, and the design can incorporate multiple candidate materials (referred to as multi-material TO). To address these TO problems, a variety of methods have been developed in literature. While notable achievements have been made, there is still potential for improvement, particularly in computational efficiency, which is primarily dictated by the need to iteratively solve the optimization problem and thereby to repeatedly solve the associated PDEs during the optimization process.

In TO, the continuous state variables (like displacement) nonlinearly depend on the binary (0/1) design variables (with 0 denoting void and 1 denoting solid), generally leading to a mixed integer nonlinear programming (MINLP) problem. One category of methods, including the Solid Isotropic Material with Penalty (SIMP) method \cite{bendsoe1999material} and the topological representative methods like Level Set \cite{wang2003level} and Phase Field \cite{takezawa2010shape}, typically relax the binary (0/1) design variables into continuous variables (varying between 0 and 1). In this way, one can avoid dealing with mixed-integer programming problems, allowing the resulting continuous nonlinear optimization problem to be solved using gradient-based nonlinear optimizers, such as the Method of Moving Asymptotes (MMA) \cite{svanberg1987method}. However, selecting an appropriate interpolation scheme to smooth the binary variable into a continuous approximation, while maintaining accuracy, can be problem-dependent. Meanwhile, introducing a nonlinear interpolation function can increase the problem's nonlinearity and add extra local minima to the resulting optimization objective  \cite{zuo2017multi, liu2024multi}. Additionally, to recover the desired binary (black/white) topology of material layout in the optimal design, additional efforts, such as Heaviside projection to convert continuous variables into binary representations as in SIMP \cite{kawamoto2011heaviside} and the floating projection (FP) constraints to enforce the grey elements evolving toward 0/1 designs as in FP \cite{huang2020smooth, huang2022three}, are required; alternatively, an additional PDE governing the evolution of the interface between void and solid phases needs to be solved, as in the level set \cite{wang2003level} and phase field \cite{takezawa2010shape} methods.

When addressing the inclusion of multiple materials, where a TO method must accurately discern the coexistence of all or some of the candidate materials in the derived material configuration, SIMP-based methods can be inefficient for the following reasons. First of all, it requires a nonlinear interpolation with respect to the normalized density and Young's modulus for each candidate material \cite{zuo2017multi}. The multiple nonlinear interpolations associated with multiple candidate materials could artificially increase the non-convexity of the resulting optimization problem by introducing extra local minima. Second, if the optimization process of a SIMP-based method starts with the inclusion of all candidate materials, the resultant material configuration will inevitably encompass them all. However, the optimal design produced from optimization should ideally exclude inefficient material(s) \cite{liu2024multi, huang2021new}, e.g., those with a lower Young's modulus or specific stiffness. As a result, a SIMP-based method must sample all possible combinations of all or some of the candidate materials, treating each as an independent TO problem. The objective function values obtained for each combination are collected and compared to determine the combination of materials with the lowest objective function value as the final optimal design. As the number of candidate materials increases, the number of TO problems that need to be solved grows exponentially, rendering this method very inefficient when dealing with multi-material designs. The FP method \cite{huang2021new} alleviates this issue by extending the FP constraints to incorporate multiple materials. This method allows solving only one TO problem with expanded design variables (by the number of candidate materials) to determine the optimal selection of materials and the corresponding topology. However, this method still requires all candidate materials sorted in a specific order according to their properties. On the other hand, when employing the level set method to solve multi-material TO problems \cite{noda2022extended}, the number of level set functions needed to represent the multi-material topology scales quadratically with the number of candidate materials. Consequently, the number of PDEs required to evolve these level set functions also scales quadratically, significantly increasing the total number of PDEs to be solved during the TO process.

The other category of methods maintain binary design variables without relaxing them into continuous variables. As TO is essentially a constrained optimization problem governed by PDEs, its solution process involves discretizing the design space, such as into finite elements, where each element represents a potential material location. Thereby, it results in an optimization problem with a large number of binary (0/1) design variables, typically one per element with a value of 1 representing the ``existence" and 0 the ``absence" of material. As the continuous state variables nonlinearly depend on the binary design variables, the TO problem in discrete setting is cast into a large-scale MINLP problem, constrained by both equality and inequality constraints. Directly solving this problem can be challenging, thus, this category of methods typically separate the binary and continuous variables into different sub-problems and then solve a sequence of these sub-problems iteratively until reaching convergence. Representative methods include Bidirectional Evolutionary Structure Optimization (BESO) \cite{huang2009bi}, Topology Optimization of Binary Structures (TOBS) \cite{picelli2021101}, Sequential Approximate Integer Programming (SAIP) \cite{liang2019topology,liang2020further}, and Generalized Benders' Decomposition (GBD) \cite{munoz2011generalized,ye2023quantum}. In BESO and TOBS, the sub-problem with respect to the binary design variables is an integer linear programming (ILP) problem. BESO employs an evolutionary algorithm to address the ILP problem \cite{huang2009bi}, but the criterion for evolving the design variables in the evolutionary algorithm depends on the physical understanding/intuition of the problem. Instead, TOBS adopts a binary optimizer to directly solve the ILP problem \cite{picelli2021101}. In SAIP, the sub-problem associated with the binary design variables is a constrained integer quadratic programming problem \cite{liang2019topology}, which is then transformed into a convex, continuous optimization problem through applying the canonical relaxation algorithm and introducing new continuous variables to replace the binary design variables. GBD was firstly introduced for truss systems \cite{munoz2011generalized} and was recently extended to address continuum TO \cite{ye2023quantum}. In GBD, the sub-problem associated with the binary design variables is a mixed integer linear programming (MILP) problem. This offers two advantages over other methods. First, it provides a deterministic optimality criterion to warrant convergence within a finite number of iterations \cite{belotti2013mixed}. Second, it enables faster convergence because the MILP is essentially a multi-cut optimization problem, where the solutions generated from all previous iteration steps are considered in the new iteration step, with each previous solution involved by one (linear) cut. In contrast, the methods like SIMP, TOBS, or SAIP only consider the solution yield from the most recent iteration. However, GBD is only applicable to convex TO problems, a condition that is rarely met in general \cite{sigmund1997design}. For dealing with multi-material TO designs, all BESO \cite{huang2009bi}, TOBS \cite{sivapuram2021design}, and SAIP \cite{liu2024multi} can be directly extended to include multiple candidate materials when solving the TO problem, whether convex or non-convex, by introducing additional sets of design variables. Each set corresponds to an additional candidate material. By doing so, the final solution of all the design variables can represent the optimal selection of materials on each discrete element. 

In light of the advantages and challenges posed by existing methods, we introduce a new TO framework in this work. This framework maintains binary design variables and addresses the large-scale MINLP problem that arises from discretizing the design space and PDE. The MINLP problem is decomposed into two sub-problems: the primal and master problems, each concerning the state variables (continuous) and the design variables (binary), respectively. By solving each sub-problem iteratively, these two types of variables are updated until convergence is achieved. In the spirit of GBD \cite{geoffrion1972generalized, munoz2011generalized, ye2023quantum}, the master problem is formulated as a multi-cut optimization problem, where each cut can be regarded as a local, linear approximation to the original objective function. To ensure the accuracy of the linear approximation within each cut, we leverage the method of trust region to constrain each cut. By doing so, the framework can tackle non-convex as well as convex TO problems, unlike GBD, which is applicable only to convex problems. Meanwhile, formulating the master problem as a multi-cut optimization problem can accelerate the convergence of the solution, because multiple linear approximations can be combined in a piecewise manner to construct a more accurate approximation of the original nonlinear objective function. Leveraging trust regions can also help address ill-conditioned optimization problems by imposing a move limit based on the change in the objective function value \cite{nocedal1999numerical}. Due to the binary nature of the design variables, the introduction of trust-region constraints only adds linear constraints, rather than quadratic ones \cite{liang2020discrete}, thereby not increasing the problem's difficulty level in solving. Furthermore, our framework allows for estimating the upper and lower bounds of the original objective function from the primal and master sub-problems, respectively. These bounds are expected to converge as the solution approaches the true value. Based on that, we establish our stopping criterion by comparing their differences during the optimization process. This way of constructing the stopping criterion is theoretically more rigorous and also practically proven to be more efficient in minimizing optimization iterations. To select appropriate trust-region radii to ensure both solution accuracy and efficiency, we allow the trust-region radius to vary adaptively during the solution process, according to the value of a merit function evaluated at the end of each iteration step. To satisfy the necessary condition for feasibility, the cuts and their corresponding trust-region constraints need to be down-selected, for which a small set of additional binary variables are introduced. This leads to a bilinear binary optimization problem, for which we propose a branching scheme to effectively handle it. Finally, to mitigate the ill-conditioning of the optimization problem resulting from extreme parameter values, we further introduce a more sophisticated parameter relaxation scheme. Different from what has been proposed in literature \cite{ye2023quantum, picelli2021101, liang2019topology, huang2010evolutionary, huang2009bi, liu2024multi}, where only one parameter is relaxed, we instead propose to relax two parameters, i.e., the minimum Young's modulus (assigned to void elements to ensure a non-singular stiffness matrix) and the target volume fraction (for single-material problems) or the target maximum permissible mass fraction (for multi-material problems). The scheme is designed such that the two parameters are relaxed in different stages and at different rates, resulting in improvements in both solution quality and efficiency.

To validate the proposed framework and demonstrate its superior performance, we present a series of numerical tests, including both the minimum compliance (convex) and compliant mechanism (non-convex) problems and for both single-material and multi-material designs. Notably, when addressing multi-material TO problems, our framework can effectively determine the optimal combination of materials without the need to either solve many independent TO problems (each corresponding to a different combination of candidate materials), sort candidate materials in a specific order, or interpolate design variables between different materials, unlike other methods \cite{zuo2017multi, huang2021new,liu2024multi}. In all cases, we assess the number of FEM analyses required to reach the optimal solution and the resulting value of the objective function, comparing these results with those obtained by other methods such as SIMP \cite{andreassen2011efficient}, FP \cite{huang2021new}, and SAIP \cite{liu2024multi}. We demonstrate that our framework can significantly reduce the number of optimization iterations, and thereby the total number of FEM analyses required throughout the solution process, while maintaining comparable solution quality (in terms of the optimal value of the objective function and the resulting material layout). Regardless of increasing discretization resolutions or the inclusion of multiple materials, both of which significantly increase the number of design variables and the latter of which also introduce more inequality constraints, we observe consistency in the solution quality and the number of optimization iterations required to reach the optimal solution. This indicates that the proposed framework is scalable with respect to the number of design variables and constraints. We therefore envision that the new TO framework proposed in this work will be especially advantageous for large-scale TO applications, which involve large numbers of design variables and constraints and require significant computational resources for FEM analyses. 

The rest of the paper is organized as follows. Section \ref{sec:problem-statement} describes the basics of the problem concerned and provides the mathematical statement of the TO problem solved in this work. In section \ref{sec:methodology}, we present in detail our proposed new framework for addressing the TO problem described in \S\ref{sec:problem-statement}. 
Section \ref{sec:numerical-results} outlines all numerical tests and their results, including single-material and multi-material minimum compliance problems and single-material and multi-material compliant mechanism designs. We compare our results to those obtained from SIMP or the state-of-the-art methods when available. The impacts of employing adaptive trust-region radii and the proposed parameter relaxation scheme are also discussed in this section. Finally, we summarize our findings and discuss potential future work in Section \ref{sec:conclusion}.
\section{Problem Statement}
\label{sec:problem-statement}

We concern topological designs of structures in continuum domains ($\Omega$) where designs are represented as distributed functions, denoted as $\rho(\mathbf{x})$ with $\mathbf{x}\in \Omega$. The design is constrained by a PDE subject to appropriate boundary conditions (BCs), which describe the underlying physics governing the design process, and inequality constraints that specify the desired material usage. Thus, the TO problem considered can be mathematically formulated as:
\begin{equation}
    \begin{aligned}
        \min_{\rho} & \quad \int_{\Omega} F[u(\mathbf{x}, \rho(\mathbf{x})), \rho(\mathbf{x})]~\mathrm{d} \Omega \\
        \text{s.t.} & \quad \mathcal{L}(u(\rho(\mathbf{x}))) + \mathbf{b} = \mathbf{0} & & \forall \mathbf{x} \in \Omega \\
        & \quad u(\mathbf{x}, \rho(\mathbf{x})) = \mathbf{u}_\Gamma(\mathbf{x}) & & \forall \mathbf{x} \in \Gamma_D \\
        & \quad \mathbf{n} \cdot \nabla u(\mathbf{x}, \rho(\mathbf{x})) = \mathbf{h}_\Gamma(\mathbf{x}) & & \forall \mathbf{x} \in \Gamma_N \\
        & \quad H_{i_H}(\rho(\mathbf{x})) \leq 0, \quad i_H = 1, \dots, n_H
    \end{aligned} \;,
    \label{eq:continuum_to}
\end{equation}
where $u$ denotes a state variable such as the displacement; $\rho = \{0,1\}$ is the design variable and describes the material distribution (either solid or void) in the design domain $\Omega$; $\mathcal{L}$ and $\mathbf{b}$ denote the differential operator and source term, respectively, in the PDE; and $\Gamma = \partial \Omega$ represents the boundary of the design domain with $\Gamma_D$ denoting the boundary where Dirichlet BC is imposed and $\Gamma_N$ the boundary where Neumann BC is enforced. The last line in Eq. \eqref{eq:continuum_to} denotes the inequality constraints exerted on the design variable $\rho$, e.g., the volume or mass constraints.

To solve the TO problem formulated in Eq. \eqref{eq:continuum_to}, numerical discretization is required, e.g., via the finite element method (FEM). In the context of designing structures governed by linear elasticity, the discretized TO problem can be written as:
\begin{equation}
    \begin{aligned}
        \min_{\mathbf{u}, \boldsymbol{\rho}} & \quad f(\mathbf{u}, \bm{\rho}) \\
        \text{s.t.} & \quad \mathbf{K}(\boldsymbol{\rho}) \mathbf{u} = \mathbf{f} \\
        & \quad H_{i_H}(\bm{\rho}) \leq 0, \quad i_H = 1, \dots, n_{H} \\
        & \quad \mathbf{u} \in \mathbb{R}^{n_u}, ~\boldsymbol{\rho} \in \{ 0, 1 \}^{n_e} 
    \end{aligned}\;,
    \label{eq:non-convex-to}
\end{equation}
where the discretized state variable $\mathbf{u} \in \mathbb{R}^{n_u}$ is evaluated on all mesh nodes, with $n_u$ the total number of mesh nodes; the discretized BCs and source term are embedded into $\mathbf{f}\in \mathbb{R}^{n_u}$; the design variable is evaluated at the centroid of each mesh element, and hence in the discrete setting, there are $n_e$ design variables, $\boldsymbol{\rho} \in \{0, 1\}^{n_e}$, to be determined from the optimization process, with $n_e$ the total number of mesh elements. With $n_u$ continuous state variables and $n_e$ binary design variables, problem \eqref{eq:non-convex-to} is a  large-scale MINLP problem. 

In Eq. \eqref{eq:non-convex-to}, the stiffness matrix is given by: 
\begin{equation}
    \mathbf{K}(\boldsymbol{\rho}) = \sum_{e=1}^{n_e} [(E_1 - E_0)\rho_e + E_0] \mathbf{K}_e \;,
    \label{eq:element_mat}
\end{equation}
where $\mathbf{K}_e \in \mathbb{R}^{{n_u} \times {n_u}}$ denotes the stiffness matrix for each mesh element under unit modulus of elasticity, with $e = 1, 2, \dots n_e$; $E_1$ is the Young's modulus of the solid material; and $E_0$ is called the \emph{minimum Young's modulus}, assigned to all void elements to guarantee that $\mathbf{K}(\boldsymbol{\rho})$ is a non-singular matrix. In practice, $E_0$ is set far smaller than the Young's modulus of the solid material, so that its impact on the solution's accuracy is negligible. In the present work, $E_0=10^{-9}$ (in reduced unit) is used, consistent with that in other methods \cite{andreassen2011efficient, huang2010evolutionary, liu2024multi}.

In this study, we particularly consider two types of inequality constraints: 1) the volume constraint to constrain the total volume of material that can be utilized within the design domain $\Omega$, which can be written as:
\begin{equation}
    \sum_{e=1}^{n_e} \dfrac{\rho_e}{n_e} \leq V_T\;,
    \label{eq:volume_constraint}
\end{equation}
with $V_T$ denoting the target total volume of material;  2) the mass constraint that limits the total mass of the designed structure, which is mainly needed when dealing with multi-material TO problems, with details provided later in \S\ref{subsec:mass_constraint}.
\section{Proposed Method}
\label{sec:methodology}
The problem formulated in Eq. \eqref{eq:non-convex-to} represents a general setting of TO problems, wherein 
the objective function may be convex or non-convex with respect to the design variables, and the inclusion of multiple candidate materials is considered. In this section, we present in detail our proposed new framework for addressing this problem. The framework is designed to significantly enhance computational efficiency, primarily by diminishing optimization iteration counts and consequently reducing the frequency of solving the associated PDE. This reduction is especially advantageous for large-scale TO applications where solving the PDE repeatedly constitutes a significant computational expense.

\subsection{Multi-cut formulation}\label{subsec:multicuts}

Inspired by the GBD method \cite{geoffrion1972generalized, munoz2011generalized, ye2023quantum}, the state variable $\mathbf{u}$ and the design variable $\bm{\rho}$ are separated into two sub-problems, i.e., the primal and master problems. By solving each sub-problem, these two variables are iteratively updated until reaching convergence. In order to derive the subproblem with respect to the design variable $\bm{\rho}$, the discretized TO problem \eqref{eq:non-convex-to} can be equivalently reformulated as:
\begin{equation}
    \begin{aligned}
        \min_{\mathbf{u}, \boldsymbol{\rho}, \bm{\mu}} & \quad f(\mathbf{u}, \bm{\rho}) + \bm{\mu}^\intercal \left(\mathbf{K}(\bm{\rho}) \mathbf{u} - \mathbf{f}\right) \\
        \text{s.t.} & \quad H_{i_H}(\bm{\rho}) \leq 0, \quad i_H = 1, \dots, n_{H} \\
        & \quad \mathbf{u} \in \mathbb{R}^{n_u}, ~\boldsymbol{\rho} \in \{ 0, 1 \}^{n_e} 
    \end{aligned}\;,
\end{equation}
where $\bm{\mu}$ is the Lagrange multiplier introduced  to move the equality constraint to the objective function. Then, in the $k$-th iteration step, the master problem determines the design variable $\bm{\rho}^k$ and is formulated as \cite{munoz2011generalized,ye2023quantum}:
\begin{equation}
    \begin{aligned}
        \min_{\boldsymbol{\rho}, \eta} & \quad \eta \\
        \text{s.t.} & \quad f(\mathbf{u}^j, \bm{\rho}^j) + \sum_{e=1}^{n_e} {(\bm{\mu}^{j})}^\intercal \mathbf{K}_e \mathbf{u}^{j} (\rho_e - \rho^{j}_e) \leq \eta, \quad j = 0, \dots, k-1 \\
        & \quad H_{i_H}(\bm{\rho}) \leq 0, \quad i_H = 1, \dots, n_{H} \\
        & \quad \boldsymbol{\rho} \in \{0, 1\}^{n_e}
    \end{aligned} \;,
    \label{eq:to_gbd_master}
\end{equation}
where, the term ${(\bm{\mu}^j)}^\intercal \mathbf{K}_e \mathbf{u}^j$ denotes the so-called sensitivity of the objective function with respect to the design variable $\bm{\rho}$ in TO. For the Lagrangian multiplier, we have $\bm{\mu} = -\mathbf{K}^{-1}(\bm{\rho}) \frac{\partial f(\mathbf{u}, \bm{\rho})}{\partial \mathbf{u}}$, which is deduced from the Karush–Kuhn–Tucker (KKT) condition \cite{nocedal1999numerical}.

In the $k$-th iteration step, the primal problem can be obtained by fixing the design variable to $\bm{\rho}^k$ and hence reads:
\begin{equation}
    \begin{aligned}
        \min_{\mathbf{u}} & \quad f(\mathbf{u}, \bm{\rho}^k) \\
        \text{s.t.} & \quad \mathbf{K}(\boldsymbol{\rho}^{k}) \mathbf{u} = \mathbf{f} \\
        & \quad \mathbf{u} \in \mathbb{R}^{n_u}
    \end{aligned} \;.
    \label{eq:to_primal}
\end{equation}
Obviously, the primal problem essentially seeks the solution of the state variable $\mathbf{u}^k$ given $\bm{\rho}^k$, i.e., 
\begin{equation}
    \mathbf{u}^k = \mathbf{K}^{-1}(\boldsymbol{\rho}^{k}) \mathbf{f} \;.
    \label{eq:to_fem}
\end{equation}
Hence, it can be directly solved by utilizing a linear system solver (direct or iterative).

Note that the master problem given in Eq. \eqref{eq:to_gbd_master} is a multi-cut optimization problem, where each cut can be regarded as a linear approximation to the original objective function. If only a single cut is included, it would be similar to the sub-problem formulated in other methods \cite{picelli2021101, liang2019topology}. The advantages of adopting the multi-cut formulation lie in the following aspects. First, when multiple cuts are combined together, we essentially construct a more accurate approximation to the original, nonlinear objective function through multiple piecewise linear approximations. Second, the solutions and sensitivity analysis obtained in multiple prior iteration steps can be leveraged for seeking the next step's solution. All these salient features expedite the solution convergence. Lastly, it also allows to employ a more rigorous convergence criterion to terminate the optimization iterations, as discussed later in \S\ref{subsec:stopping_criteria}.

To ensure the accuracy of the linear approximation within each cut, we require the expected solution in each cut not too far from $\bm{\rho}^j$, i.e., $(\rho_e - \rho^{j}_e)$ in Eq. \eqref{eq:to_gbd_master} needs to be constrained. For that, we leverage the concept of trust region, as discussed in the next section.

\subsection{Trust region}
\label{subsec:trust_region}

The trust region for a given point $\bm{\rho}^j$ can be defined with a radius $d^j$ and constructed through the $L_2$ norm, following \cite{liang2020discrete}, as:
\begin{equation}
    \dfrac1{n_e} \sum_{e = 1}^{n_e} | \rho_e - \rho_e^j |^2 \leq d^j \;.
    \label{eq:trust_region}
\end{equation}
One can observe that
\begin{equation}
| \rho_e - \rho_e^j |^2 = \rho_e^2 - 2\rho_e \rho_e^j + {\rho_e^j}^2 = (1 - 2\rho_e^j) \rho_e + {\rho_e^j}^2 \;,
\end{equation}
resulting from the fact that $\bm{\rho}$ is binary. Eq. \eqref{eq:trust_region} can then be rewritten as a linear constraint given by:
\begin{equation}
    t^j(\bm{\rho}) = \dfrac1{n_e} \left[\sum_{e = 1}^{n_e} \sigma(\rho_e^j) \rho_e + \sum_{e = 1}^{n_e} {\rho_e^j}^2 \right] \leq d^j \;,  \quad \quad \text{with} \quad \sigma(\rho_e^j) = 1 - 2 \rho_e^j \;.
    \label{eq:move_limit}
\end{equation}
By incorporating the trust-region constraint given by Eq. \eqref{eq:move_limit}, the master problem in Eq. \eqref{eq:to_gbd_master} can be reformulated as: 
\begin{equation}
    \begin{aligned}
        \min_{\boldsymbol{\rho}, \eta} & \quad \eta \\
        \text{s.t.} & \quad \widetilde{f}^j(\bm{\rho}) \leq \eta, \quad j = 0, \dots, k - 1 \\
        & \quad t^j(\bm{\rho}) \leq d^j, \quad j = 0, \dots, k - 1 \\
        & \quad H_{i_H}(\bm{\rho}) \leq 0, \quad i_H = 1, \dots, n_H \\
        & \quad \boldsymbol{\rho} \in \{0, 1\}^{n_e} 
    \end{aligned}\;,
    \label{eq:multicut_trust_region}
\end{equation}
where we define $\widetilde{f}^j(\bm{\rho})= f(\mathbf{u}^j, \bm{\rho}^j) + \sum_{e=1}^{n_e} {(\bm{\mu}^{j})}^\intercal \mathbf{K}_e \mathbf{u}^{j} (\rho_e - \rho^{j}_e)$ for convenience of notation.

The problem stated in Eq. \eqref{eq:multicut_trust_region} has a solution only if all trust regions share a common region, which is a necessary condition for feasibility. To satisfy this condition, the multiple cuts included in Eq. \eqref{eq:multicut_trust_region} need to be down-selected. For that, we introduce an additional set of binary variables, denoted as $\bm{\alpha} = \{ \alpha_1, \alpha_2, \cdots, \alpha_k \}$, and reformulate Eq. \eqref{eq:multicut_trust_region} as:
\begin{equation}
    \begin{aligned}
        \min_{\boldsymbol{\rho}, \bm{\alpha}, \eta} & \quad \eta \\
        \text{s.t.} & \quad \widetilde{f}^j(\bm{\rho}) \alpha_j - \Pi(1 - \alpha_j) \leq \eta, \quad j = 0, \dots, k-1 \\
        & \quad t^j(\bm{\rho}) \alpha_j \leq d^j, \quad j = 0, \dots, k-1 \\
        & \quad H_{i_H} (\bm{\rho}) \leq 0, \quad i_H = 1, \dots, n_H \\
        & \quad \sum_{j=0}^{k-1} \alpha_j \geq 1 \\
        & \quad \bm{\alpha} \neq \bm{\alpha}^i, \quad i = 1, \dots, k - 1 \\
        & \quad \bm{\rho} \in \{0, 1\}^{n_e}, \bm{\alpha} \in \{0, 1\}^k \;.
    \end{aligned}
    \label{eq:multicut_trustregion_activecut}
\end{equation}
where $\Pi$ is a large positive real number; only the cuts with $\alpha_j = 1$ are effective and hence are called ``active" cuts. Here, $\bm{\alpha}^i$ with $i = 1, \dots, k - 1$ denotes the solutions of $\bm{\alpha}$ to problem \eqref{eq:multicut_trustregion_activecut} obtained from previous iteration steps. By constraining $\bm{\alpha} \neq \bm{\alpha}^i$, all previous solutions of $\bm{\alpha}$ are excluded in the current $k$-th iteration step, so that the optimization for $\bm{\rho}$ can be advanced rather than stagnant at a previously attained solution.

Furthermore, to mitigate the so-called ``checkerboard" artifact arising from finite element analysis (FEA) with uniform meshing, filtering is needed, which can be realized by modifying the objective function in the master problem \eqref{eq:multicut_trustregion_activecut} as:
\begin{equation}
    \widetilde{f}^j (\bm{\rho}) = f(\mathbf{u}^j, \bm{\rho}^j) + \sum_{e=1}^{n_e} {\widetilde{w}_e^j} (\rho_e - \rho^{j}_e) \;,
    \label{eq:linear_approximation}
\end{equation}
where
\begin{equation}
    \widetilde{w}_e^j = \dfrac{\sum_{e^\prime \in \mathcal{N}_e^r} h_{e, e^\prime}(r) w_e^j}{\sum_{e^\prime \in \mathcal{N}_e^r} h_{e, e^\prime}(r)}\;,
    \label{eq:sensitivity_filtering}
\end{equation}
with 
\begin{equation}
    w_e^j = [(E_1 - E_0) \rho_e^j + E_0] (\bm{\mu}^j)^\intercal \mathbf{K}_e \mathbf{u}^j
    \label{eq:sensitivity-regularization}
\end{equation}
and $h_{e, e^\prime}(r) = \max(0, r - \|\mathbf{x}_e - \mathbf{x}_{e^\prime}\|_2)$, following \cite{picelli2021101,sun2022sensitivity}. 

\subsection{Algorithm for solving the master problem \eqref{eq:multicut_trustregion_activecut}}
\label{subsec:sovling_bilinear}

The master problem as described in Eq. \eqref{eq:multicut_trustregion_activecut} contains the nontrivial bilinear terms formed by the binary variables $\bm{\rho}$ and $\bm{\alpha}$. As a result, it cannot be directly solved by employing an off-the-shelf integer programming solver. Thus, we herein propose a feasible and efficient way to solve it. By noting that the size of $\bm{\alpha}$ is much smaller than that of $\bm{\rho}$, branching on $\bm{\alpha}$ is not difficult. The branching essentially selects one or multiple indices $j$, $j \in \{0, 1, \cdots, k-1\}$, such that $\alpha_j=1$ but the rest elements in $\bm{\alpha}$ are zero. For example, 1 and 2 are selected, leading to $\bm{\alpha}=(0,1,1,0,\dots)$. By doing so, the constraint $\sum_{j = 0}^{k-1} \alpha_j \geq 1$ is satisfied. Further, the selection must also satisfy the constraint: $\bm{\alpha} \neq \bm{\alpha}^i$, $i = 1, 2, \cdots, k - 1$. Let $\mathcal{P}_s \subseteq \{ 0, 1, \cdots, k - 1 \}$ denote one selection for the index $j$ that satisfy both constraints, and all such selections are collected into the set $\mathcal{C}(k)$.

Among all selections in $\mathcal{C}(k)$, there is only one selection (denoted as $\mathcal{P}_1$), for which only one element of $\bm{\alpha}$ is nonzero, and this only nonzero element is $\alpha_{k-1}$, i.e., $\mathcal{P}_1 =\{k - 1 \}$ and $\bm{\alpha}=(0,0,\dots,1)$. (Note that $\bm{\alpha}$ with a single nonzero element $\alpha_j$ for $\forall j<k-1$ has been considered in previous iterations and hence would not be considered again in the current $k$-th iteration step.) For this particular selection, the master problem in Eq. \eqref{eq:multicut_trustregion_activecut} is reduced to a single-cut problem and can be simplified as: 
\begin{equation}
    \begin{aligned}
        \min_{\boldsymbol{\rho}, \eta} & \quad \eta \\
        \text{s.t.} & \quad \widetilde{f}^{k-1}(\bm{\rho}) \leq \eta \\
        & \quad t^{k-1}(\bm{\rho}) \leq d^{k-1} \\
        & \quad H_{i_H} (\bm{\rho}) \leq 0, \quad i_H = 1, \dots, n_H \\
        & \quad \boldsymbol{\rho} \in \{0, 1\}^{n_e} \;.
    \end{aligned}
    \label{eq:single-cut}
\end{equation}
We denote its solution as $\eta_{k, 1}$. (The solution for each of the other single-cut problems that have been solved in the previous iteration steps is denoted as $\eta_{l, 1}$ with $\forall l<k$.) 

For every other selection in $\mathcal{C}(k)$, $|\mathcal{P}_s|\ge 2$, i.e., $\bm{\alpha}$ has at least two nonzero elements. For that, the master problem in Eq. \eqref{eq:multicut_trustregion_activecut} must involve multiple cuts (at least two) and can be rewritten as:
\begin{equation}
    \begin{aligned}
        \min_{\boldsymbol{\rho}, \eta} & \quad \eta \\
        \text{s.t.} & \quad \widetilde{f}^j(\bm{\rho}) \leq \eta, \quad \forall j \in \mathcal{P}_s \\
        & \quad t^j(\bm{\rho}) \leq d^j, \quad \forall j \in \mathcal{P}_s \\
        & \quad H_{i_H}(\bm{\rho}) \leq 0, \quad i_H = 1, \dots, n_H \\
        & \quad \boldsymbol{\rho} \in \{0, 1\}^{n_e} \;.
    \end{aligned}
    \label{eq:branched_multi-cut}
\end{equation}
We denote its solution as $\eta_{k, s}$. As the multi-cut problem includes more than one cuts (or constraints), its solution $\eta_{k, s}$ cannot be lower than the solution of any single-cut problem containing just one of those cuts, i.e., $\eta_{k, s} \ge \max_{l \in \mathcal{P}_s} (\eta_{l, 1})$. We can utilize this feature for early stopping the branching on $\bm{\alpha}$, thereby enhancing computational efficiency. 

To proceed, we denote $\bar{\eta}_{k,s} = \max_{l \in \mathcal{P}_s} (\eta_{l, 1})$. All selections represented by $\mathcal{P}_s$ with $|\mathcal{P}_s|\ge 2$ in $\mathcal{C}(k)$ are ranked in non-descending order based on the value of $\bar{\eta}_{k,s}$. Starting from the first in the rank list, we solve the corresponding multi-cut problem, as formulated in Eq. \eqref{eq:branched_multi-cut}, until we find that the minimum of the solutions of all solved multi-cut problems is smaller than $\bar{\eta}_{k,s}$ of the next multi-cut problem in the rank list, i.e., $\min_s(\eta_{k, s}) < \bar{\eta}_{k,s+1}$. Once this is achieved, branching on $\bm{\alpha}$ can be terminated, because the remaining selections in the list would not yield better solutions. Thus, the final solution of $\eta$ for the master problem \eqref{eq:multicut_trustregion_activecut} is given by:
\begin{equation}
    \eta^k = \min\left[\min_s(\eta_{k, s}),~\eta_{k, 1}\right] \;,
\end{equation}
i.e., the minimum among the solutions of all solved multi-cut problems (as in Eq. \eqref{eq:branched_multi-cut}) and the single-cut problem (as in Eq. \eqref{eq:single-cut}). The final solution of $\bm{\alpha}$ and $\boldsymbol{\rho}$ will then be derived from the multi-cut or single-cut problem that yields the minimum $\eta$.

If the final solution for the master problem \eqref{eq:multicut_trustregion_activecut} is determined as the solution of a multi-cut problem, i.e., $\eta^k = \min_s(\eta_{k, s})$, it indicates that using a single cut solely based on the solution of the last iteration step, as done by the methods like SIMP \cite{andreassen2011efficient} and TOBS \cite{picelli2021101}, does not always yield the best solution; but including multiple cuts based on the solutions from at least two previous iteration steps can improve the solution. Therefore, allowing for adaptively including multiple cuts would enable us to find the best solution at each iteration step, thereby accelerating convergence and improving solution efficiency. This highlights the distinction of the framework proposed in this paper from other methods, including SIMP \cite{andreassen2011efficient}, TOBS \cite{picelli2021101}, FP \cite{huang2020smooth, huang2021new}, and DVTOPCRA \cite{liang2019topology, liu2024multi}.

In summary, solving the master problem as formulated in Eq. \eqref{eq:multicut_trustregion_activecut} involves solving two types of linear integer programming problems, i.e., the single-cut problem in Eq. \eqref{eq:single-cut} and the multi-cut problems in Eq.  \eqref{eq:branched_multi-cut}. Each type of problem can be directly solved using a linear integer programming solver; in this work, we use the Gurobi optimizer \cite{gurobi}. The algorithm proposed for solving the master problem \eqref{eq:multicut_trustregion_activecut} is summarized in \Cref{algorithm:multi_cuts_solver}. 
\begin{algorithm}
    \caption{\textsc{MasterProbSolver}($\bm{\alpha}^i$)}
    \label{algorithm:multi_cuts_solver}
    \textbf{Input}: $\bm{\alpha}^i$ from previous iteration steps, $i = 1, \dots, k - 1$ 
    
    \textbf{Output}: The optimal solution $\bm{\alpha}^k$, the minimizer $\bm{\rho}^k$ and the minimum value $\eta^k$
    
    \textbf{Utility}: Solve the master problem \eqref{eq:multicut_trustregion_activecut}
    
    \begin{algorithmic}[1]
        \State Form the set $\mathcal{C}(k)$ for branching on $\bm{\alpha}$
        \State Solve the single-cut problem in Eq. \eqref{eq:single-cut} and yield the solution $\bm{\rho}_{k, 1}$ and  $\eta_{k, 1}$
        \State Sort $\mathcal{P}_s\in \mathcal{C}(k)$ (with $|\mathcal{P}_s|\ge 2$) in non-descending order based on the value of $\bar{\eta}_{k,s}$ 
        \For {$s = 2, \dots, |\mathcal{C}(k)|$}
            \State Solve the multi-cut problem in Eq. \eqref{eq:branched_multi-cut} and yield the solution $\bm{\rho}_{k, s}$ and  $\eta_{k, s}$ 
            \If {$\min_s(\eta_{k, s}) < \bar{\eta}_{k,s+1}$}
                \Break
            \EndIf
        \EndFor
        \State $s^* = \argmin_{s} (\eta_{k,s})$
        \State $\boldsymbol{\rho}^k=\bm{\rho}_{k, s^*}$, $\eta^k = \eta_{k,s^*}$
        \State Let $\alpha^k_l = 1,~\forall l \in \mathcal{P}_{s^*}$\;; ~~$\alpha^k_l = 0,~\forall l \notin \mathcal{P}_{s^*}$
    \end{algorithmic}
    \hspace*{\algorithmicindent} \textbf{Return}: $\bm{\alpha}^k$, $\boldsymbol{\rho}^k$,  $\eta^k$
\end{algorithm}

\subsection{Adaptive trust-region radius}
\label{subsec:adaptive_adjustment}

Employing appropriate trust-region radii in Eq. \eqref{eq:move_limit} is crucial for ensuring both solution accuracy and efficiency. Thus, it can be advantageous to vary it adaptively during the solution process when solving the master problem \eqref{eq:multicut_trustregion_activecut}. By leveraging a merit function as given by \cite{yuan2015recent}: 
\begin{equation}
    \omega^k =  \min_{j \in \mathcal{P}_{s^*}(k)} \left(\dfrac{f(\mathbf{u}^{j}, \bm{\rho}^{j}) - f(\mathbf{u}^{k}, \bm{\rho}^{k})}{f(\mathbf{u}^{j}, \bm{\rho}^{j}) - \eta^k}\right) \;,
    \label{eq:merit_fun}
\end{equation}
the trust-region radius $d^k$ can be updated at each iteration step $k$, following the rule:
\begin{equation}
    d^k =
    \left\{
        \begin{aligned}
            & \max(\theta_1 d^*, d_{\min}), &\quad &0 \leq \omega^k < 1 \\
            & \min(\theta_2 d^*, d_{\max}), &\quad &\omega^k \geq 1 \\
            & \max(0.5 d^*, d_{\min}), &\quad &\omega^k < 0 \\
        \end{aligned}
    \right.
    \label{eq:trust_region_factor}
\end{equation}
with $d^* = \min_{j \in \mathcal{P}_{s^*}(k)} (d^j)$. Here, $\theta_1$ and $\theta_2$ are the adjustment (shrinking/enlarging) factors for the trust-region radius and are assigned as 0.7 and 1.5, respectively, in all numerical tests. And we adopt $d_{\min} = 10^{-3}$ (sufficiently small) and $d_{\max} = 0.6$ (sufficiently large) in this paper.

The merit function defined in Eq. \eqref{eq:merit_fun} evaluates the ratio between the exact reduction in the objective function and the reduction estimated from problem \eqref{eq:multicut_trustregion_activecut}. The desired trust-region radius $d^k$ renders the merit function satisfying $\omega^k \geq 1$. Thus, if $\omega^k \geq 1$, $d^k$ may be enlarged, for which a magnification factor 1.5 is employed to enlarge the trust-region radius for the newly added $k$-th cut for the next iteration step. Conversely, if $0 \leq \omega^k < 1$, the exact reduction in the objective function is smaller than the one estimated from problem \eqref{eq:multicut_trustregion_activecut}, and hence $d^k$ needs to be decreased. Here, a shrinking factor 0.7 is employed to reduce the trust-region radius for the $k$-th cut. However, if $\omega^k < 0$, it implies that the radius is too large that the searching goes to the wrong direction. Hence, a more aggressive shrinkage for the trust-region radius is demanded, and here we empirically choose a factor of 0.5.

\subsection{Algorithm for the entire optimization process}\label{subsec:stopping_criteria}

The iterations of solving the primal problem \eqref{eq:to_primal} and the master problem \eqref{eq:multicut_trustregion_activecut} yield a series of upper bounds and lower bounds, respectively, for the objective function in the original problem \eqref{eq:non-convex-to}. The gap between the lowest upper bound $U = \min_{j}(f(\mathbf{u}^j, \bm{\rho}^j))$ with $j\leq k$ and the lower bound $\eta^k$ should be minimized until reaching  $\tfrac{|\eta^k - U|}{|U|} < \varepsilon$ (with $\varepsilon$ the preset tolerance) to terminate the optimization iterations. Additionally, the criterion $\eta^k > U$ is enforced, meaning the iterations should stop once the lower bound exceeds the upper bound, i.e., when the master problem \eqref{eq:multicut_trustregion_activecut} finds a solution greater than $U$. The entire framework proposed for solving the TO problem in Eq. \eqref{eq:non-convex-to} is summarized in Algorithm \ref{algo:multi-cuts}.
\begin{algorithm}[htp]
    \caption{\textsc{TOSolver}($\bm{\rho}^0$)}
    \label{algo:multi-cuts}
    \textbf{Input}: Initial material configuration $\bm{\rho}^0$
    
    \textbf{Output}: Optimal material layout $\bm{\rho}^*$
    
    \begin{algorithmic}[1]
        \State Employ a linear system solver to obtain $\mathbf{u}^{0} = \mathbf{K}^{-1}(\boldsymbol{\rho}^{0}) \mathbf{f}$
        \State Evaluate the objective function value $f(\mathbf{u}^0, \bm{\rho}^0)$
        \State $\bm{\rho}^* = \bm{\rho}^0$, $U = \infty$
        \For {$k=1,\cdots$}
            \State $\bm{\alpha}^k, \bm{\rho}^k, \eta^k \gets \textsc{MasterProbSolver}(\bm{\alpha}^i)$ 
            \State Employ a linear system solver to obtain $\mathbf{u}^{k} = \mathbf{K}^{-1}(\boldsymbol{\rho}^{k}) \mathbf{f}$
            \State Evaluate the objective function value $f(\mathbf{u}^{k}, \bm{\rho}^k)$
            \If { $f(\mathbf{u}^k, \bm{\rho}^k) < U$ }
                \State $U = f(\mathbf{u}^k, \bm{\rho}^k)$, $\bm{\rho}^* = \bm{\rho}^k$
            \EndIf
            \If { $\tfrac{|\eta^k - U|}{|U|} < \varepsilon$ or $\eta^k > U$ }
                \Break
            \EndIf
            \State Adjust the trust region radius according to Eq. \eqref{eq:trust_region_factor}
        \EndFor
    \end{algorithmic}
    \hspace*{\algorithmicindent} \textbf{Return}: $\boldsymbol{\rho}^*$
\end{algorithm}

It is worth noting that the stopping criterion used here differs from those in other methods. In SIMP \cite{andreassen2011efficient}, the optimization iterations are terminated when the changes in the design variables become small over multiple consecutive iterations. In the methods that directly handle binary design variables and draw on integer optimization \cite{huang2009bi, picelli2021101, liang2019topology}, including the framework proposed in this work, the stopping criterion is typically based on the objective function values. However, BESO \cite{huang2009bi}, TOBS \cite{picelli2021101}, and SAIP \cite{liang2019topology} terminate the optimization iterations when the changes in the objective function values become small over multiple consecutive iterations, usually requiring at least 10 consecutive iterations. In contrast, our framework allows for the estimation of the upper and lower bounds of the original objective function by decomposing the problem into primal and master sub-problems, which are expected to converge as the solution nears the true value. Thus, we establish our stopping criterion by comparing the difference between these upper and lower bounds, which is both theoretically sounder and practically more efficient.

\subsection{Extension to multi-material TO}
\label{subsec:mass_constraint}
For TO designs aiming to incorporate multiple candidate materials, the design variable's space expands, i.e., $\bm{\rho} \in \{0, 1\}^{n_e\times n_M}$, where $n_M$ is the total number of material candidates considered. The objective function in Eq. \eqref{eq:multicut_trustregion_activecut} then needs to be extended as:
\begin{equation}
    \widetilde{f}^j (\bm{\rho}) = f(\mathbf{u}^j, \bm{\rho}^j) + \sum_{e=1}^{n_e} \left( \sum_{m=1}^{n_M} \left[ {\widetilde{w}_{e,m}^j} (\rho_{e, m} - \rho^{j}_{e,m}) \right] \right) \;,
    \label{eq:linear_approximation_multimat}
\end{equation}
where the index $m = 1, 2, \dots, n_M$ corresponds to each candidate material.

The mass constraint that limits the total mass of the designed structure, following the discrete multi-material optimization (DMO) scheme \cite{stegmann2005discrete}, can be expressed as:
\begin{equation}
    \begin{aligned}
        & \sum_{e=1}^{n_e} \sum_{m=1}^{n_M} \frac{\bar{M}_m}{n_e} \rho_{e, m} \leq \bar{M}_{\max} \\
        & \sum_{m=1}^{n_M} \rho_{e, m} \leq 1, \quad e = 1, 2, \dots, n_e 
    \end{aligned}\;,
    \label{eq:mass_constraint}
\end{equation}
where $\bar{M}_{\max}$ is the target maximum total mass, and $\bar{M}_m$ is the normalized density for each candidate material. This mass constraint consists of two components: the first one imposes that the total mass of the resulting structure can be no greater than the allowed maximum mass $\bar{M}_{\max}$; the second enforces that each mesh element $e$ can contain at most one candidate material or be void, but not multiple materials simultaneously. In total, the mass constraint specified in Eq. \eqref{eq:mass_constraint} encompasses $1+n_e$ inequality constraints. Note that the mass constraint can be regarded as a generalization of the volume constraint: if only a single material is involved, and the normalized density is $\bar{M}_1 = 1.0$, the first line of the above mass constraint is equivalent to the volume constraint in Eq. \eqref{eq:volume_constraint} when setting $\bar{M}_{\max} = V_T$, and the second line is automatically satisfied and hence not needed. It is also worth noting that when more than one candidate material is involved, the design variables increase by a factor of $n_M$ (resulting in a total of $n_e \times n_M$ design variables), and the number of constraints also scales with the number of design variables, as the second line in Eq. \eqref{eq:mass_constraint} essentially consists of $n_e$ inequality constraints. Therefore, the greatly expanded design space and increased number of inequality constraints associated with the problem make multi-material TO design a challenging task.

Furthermore, the stiffness matrix, sensitivity, and trust regions also need to be adjusted accordingly. However, unlike other methods \cite{zuo2017multi, huang2021new,liu2024multi}, our approach eliminates the necessity for sorting candidate materials in a specific order and avoids any interpolation of design variables between different materials. Without these additional efforts, we modify the stiffness matrix in Eq. \eqref{eq:element_mat} as:
\begin{equation}
    \mathbf{K}(\bm{\rho}) = \sum_{e=1}^{n_e} \sum_{m=1}^{n_M} \left[(E_m - E_0) \rho_{e, m} + E_0 \right] \mathbf{K}_{e} \;.
    \label{eq:multi-material_stiffness}
\end{equation}
And the sensitivity is given by:
\begin{equation}
    w^j_{e, m} =
    \left\{
    \begin{aligned}
        & (E_m - E_0) (\bm{\mu}^j)^\intercal \mathbf{K}_e \mathbf{u}^j, & \quad & \text{ if } \rho_{e, m}^j = 1 \\
        & E_m (E_{m^\prime} - E_0) (\bm{\mu}^j)^\intercal \mathbf{K}_e \mathbf{u}^j, & \quad & \text{ if } \rho_{e, m}^j = 0, ~\rho_{e, m^\prime}^j=1, ~ m'\neq m \\
        & E_0 (\bm{\mu}^j)^\intercal \mathbf{K}_e \mathbf{u}^j, & \quad & \text{ if } \sum_{m^\prime = 1}^{n_M} \rho_{e, m^\prime}^j = 0
    \end{aligned}\;,
    \right.
    \label{eq:sensitivity_filtering_multimaterial}
\end{equation}
where $m'$, another index for candidate materials, is used to distinct from $m$. Eq. \eqref{eq:sensitivity_filtering_multimaterial} can be regarded as a generalization of the sensitivity evaluation in Eq. \eqref{eq:sensitivity-regularization} to multi-material TO. Here, the specification for sensitivity in Eq. \eqref{eq:sensitivity_filtering_multimaterial} implicitly embeds material preferences in the selection process. The first line specifies the sensitivity for the element $e$ to keep the same candidate material. The second line determines the sensitivity for the element $e$ to switch from candidate material $m'$ to $m$. Assigning the same sensitivity for a void element transitioning into a solid element with any material (third line of Eq. \eqref{eq:sensitivity_filtering_multimaterial}) implies that the material with the smallest normalized density will be selected first during optimization. Finally, filtering is still necessary to mitigate the ``checkerboard" artifact. Thus, the filtered sensitivity $\widetilde{w}^j_{e, m}$ for each candidate material $m$ can be obtained in the same manner as described in Eq. \eqref{eq:sensitivity_filtering}.

The trust-region construction can be extended as:
\begin{equation}
    \frac{1}{n_M n_e} \sum_{e = 1}^{n_e} \left| \sum_{m=1}^{n_M} \rho_{e, m} - \sum_{m=1}^{n_M} \rho_{e, m}^j \right|^2 \leq d^j \;.
    \label{eq:multi-material-trust_region}
\end{equation}
Since any given element can take at most one solid material (as also stated by the multi-material mass constraint in Eq. \eqref{eq:mass_constraint}), the bilinear term $\rho_{e, m} \rho_{e, m^\prime}$ with $m\neq m'$ derived from $\left( \sum_{m=1}^{n_M} \rho_{e, m} \right)^2$ is always zero. As a result, the trust region extended to multi-material scenarios remains a linear constraint with respect to the design variable $\rho_{e, m}$, as can be seen from:
\begin{equation*}
     \sum_{e = 1}^{n_e} \left| \sum_{m=1}^{n_M} \rho_{e, m} - \sum_{m=1}^{n_M} \rho_{e, m}^j \right|^2 = \sum_{e = 1}^{n_e} \left( 1 - 2 \sum_{m=1}^{n_M} \rho_{e, m}^j \right) \left( \sum_{m=1}^{n_M} \rho_{e, m} \right) + \sum_{e = 1}^{n_e} \left( \sum_{m=1}^{n_M} \rho_{e, m}^j \right)^2 \;.
\end{equation*}

\subsection{Parameter relaxation scheme}
\label{subsec:parameter_relaxation}
Due to the small value of the minimum Young's modulus assigned to void elements (i.e., $E_0 = 10^{-9}$), and the small target volume fraction (e.g., $V_T = 0.3$) or small maximum permissible total mass fraction (e.g., $\bar{M}_{\max}$ = 0.3), the optimization problem can be ill-conditioned, as evidenced in \S\ref{subsubsec:parameter-adjustment-improvement}. Thus, different TO methods usually relax some parameter's values during the solution process, i.e., starting with a relatively larger value and later reducing the parameter's value to the desired one, in order to ease the ill-conditioning issue. Typically, for single-material problems, the target volume fraction is chosen to relax, e.g., in \cite{ye2023quantum, picelli2021101, liang2019topology, huang2010evolutionary}; for multi-material problems, the target maximum permissible mass fraction is typically relaxed, e.g., in \cite{huang2009bi, liu2024multi}.

In the present work, we further explore relaxing the minimum Young's modulus $E_0$, alongside the target volume fraction $V_T$ or the target maximum permissible mass fraction $\bar{M}_{\max}$. From numerical experiments, we find that the minimum Young's modulus contributes the most to the conditioning of the optimization problem, with more details provided in \S \ref{subsubsec:parameter-adjustment-improvement}. Thus, we consider the minimum Young's modulus as the primary parameter for relaxation and the target volume fraction or the maximum permissible total mass fraction as the secondary parameter for relaxation. Accordingly, the proposed parameter relaxation scheme consists of two parts. In the first part, the value of the target volume fraction or the maximum permissible total mass fraction is changed, while the minimum Young's modulus is fixed at a larger value $10^{-2}$. This part can consist of $N_P$ stages. In each stage, the parameter's value is gradually changed from an initial larger value towards the desired one, following:
\begin{equation}
    P_l = P_0 e^{-\frac{l - 1}{N_P - 1} \log \frac{P_0}{P_d}}, \quad l = 1, \dots N_P \;,
    \label{eq:relaxation_interpolation}
\end{equation}
where $P_0$ denotes the initial value for the parameter; $P_d$ denotes the desired value of the parameter; and $P$ can be either the target volume fraction $V_T$ or the maximum permissible total mass fraction $\bar{M}_{\max}$. Here, an exponential function is employed to achieve a substantial initial reduction in the parameter values, followed by a gradual and smoother decrease, ultimately reaching the desired value. The second part consists of only one stage, where $E_0$ is varied from $10^{-2}$ to $10^{-9}$, while the value of the target volume fraction or the maximum permissible total mass fraction is fixed at its desired value. As a result, the entire parameter relaxation scheme involves $N_P + 1$ stages in total. Such proposed parameter relaxation scheme is summarized in Algorithm \ref{algo:parameter_reduction}.
\begin{algorithm}[htp]
    \caption{\textsc{ParameterRelaxation}($\mathbf{X}$)}
    \label{algo:parameter_reduction}
    \textbf{Input}: $\mathbf{X}$, which are the minimum Young's modulus $\mathbf{E} = \{E_0^1, \cdots, E_0^{N_P+1}\}$ and the target volume fraction $\mathbf{V} = \{ V_1, \cdots, V_{N_P+1} \}$ or the maximum permissible total mass fraction $\bar{\mathbf{M}} = \{ \bar{M}_1, \cdots, \bar{M}_{N_P+1} \}$
    
    \textbf{Output}: Optimal material layout $\bm{\rho}^*$
    
    \begin{algorithmic}[1]
        \State Initialize the design variables as $\bm{\rho}^{*} = V_T \mathbf{1}$ or $\bm{\rho}^{*} = \bar{M}_{\max} \mathbf{1}$, where $\mathbf{1}$ is a unit vector
        \State Initialize the trust-region radius $d^0$
        \For {$l=1,\cdots N_P+1$}
        \State $\boldsymbol{\rho}^{0} \gets \boldsymbol{\rho}^*$
            \If {$l\leq N_P$}
                \State $E_0 \gets E_0^l=10^{-2}$ 
                \State Evaluate $V_l$ or $\bar{M}_l$ from Eq. \eqref{eq:relaxation_interpolation}
                \State $V_T \gets V_l$ or $\bar{M}_{\max} \gets \bar{M}_l$ 
            \EndIf
            \If{$l = N_P+1$}
                \State $E_0 \gets E_0^l=10^{-9}$ 
                \State $V_T \gets V_l=V_d$ or $\bar{M}_{\max} \gets \bar{M}_l=\bar{M}_d$
            \EndIf
            \State $\bm{\rho}^* \gets \textsc{TOSolver}(\bm{\rho}^0)$
        \EndFor
    \end{algorithmic}
    \hspace*{\algorithmicindent} \textbf{Return}: $\boldsymbol{\rho}^*$
\end{algorithm}
In the present work, for the very first stage, the initial material configuration $\bm{\rho}^0$ is set to a gray material layout, same as that in SIMP. For subsequent stages, the initial material configuration $\bm{\rho}^0$ in Algorithm \ref{algo:parameter_reduction} is the resulting binary (black/white) material layout yield at the end of the previous stage.

\section{Numerical Results}
\label{sec:numerical-results}
In this section, we systematically assess the new framework introduced in \S \ref{sec:methodology} by solving different TO problems, focusing on two primary types: minimum compliance problems and compliant mechanism problems. For each type of problems, we examined a variety of scenarios involving single-material or multi-material designs and different design domains, BCs, discretization resolutions, and material properties. The results are compared with those obtained by employing the conventional SIMP method and the methods lately reported in literature, including SAIP with recursive multiphase materials interpolation \cite{liu2024multi} and FP \cite{huang2021new}. The comparison concerns both the number of FEM analyses required during the solution process and the optimal value achieved for the objective function. By doing so, we demonstrate in a systematic way the superior performance of our proposed new TO framework.

\subsection{Single-material Minimum Compliance}
\label{subsec:single-mc}
\begin{figure*}[ht]
    \centering
    \begin{subfigure}[t]{.45\textwidth}
        \centering
        \includegraphics[width=\textwidth]{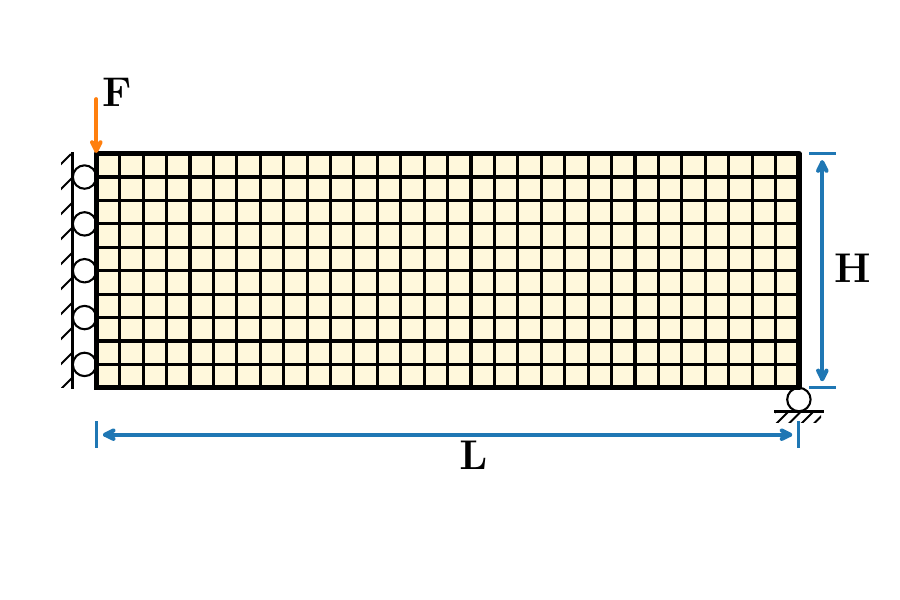}
        
        \vspace{-0.2in}
        \caption{MBB}
        \label{subfig:mbb}
    \end{subfigure}
    \quad
    \begin{subfigure}[t]{.45\textwidth}
        \centering
        \includegraphics[width=\textwidth]{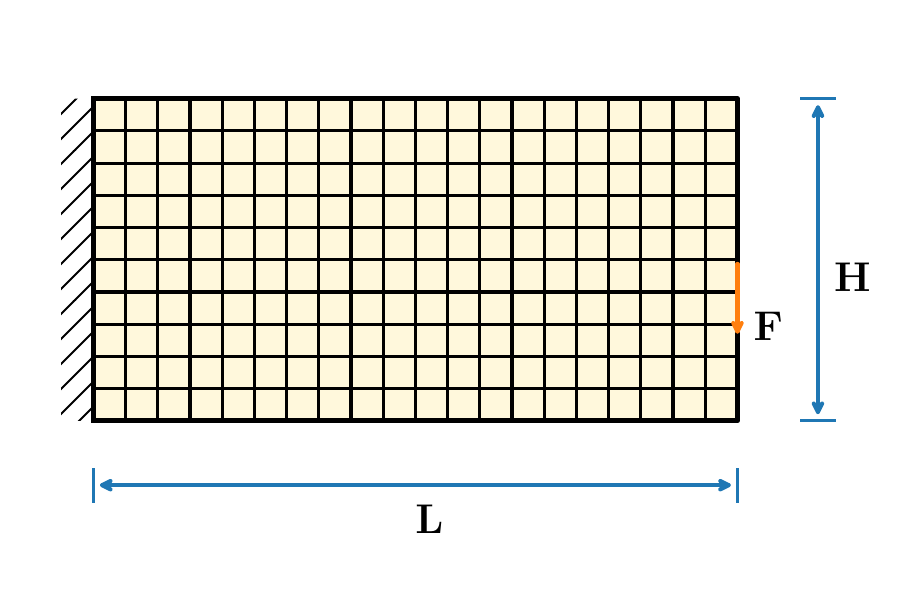}
        
        \vspace{-0.2in}
        \caption{Cantilever}
        \label{subfig:cantilever}
    \end{subfigure}
    \caption{Single-material minimum compliance: design domain and boundary conditions for two different problems.}
    \label{fig:minimum_compliance}
\end{figure*}

We first consider the minimum compliance design for a single material, but with two different settings, as depicted in \Cref{fig:minimum_compliance}. The first setting corresponds to the Messerschmitt-Bölkow-Blohm (MBB) beam problem, as illustrated in \Cref{subfig:mbb}. It is designated in a rectangular domain with a length-height ratio of $L:H=3:1$. The left-side edge of the rectangular domain is constrained in $x$ (horizontal) direction but can freely move along $y$ (vertical) direction. The movement of the bottom-right corner is restricted solely along $x$-axis. An external force $\mathbf{F}_y = -1$ is exerted at the top-left corner of the design domain. The second setting represents the design of a cantilever, defined in a rectangular domain with a length-height ratio of $L:H=2:1$, as shown in \Cref{subfig:cantilever}. The left side of the rectangular domain is constrained in both $x$ and $y$ axes. In both problems, the goal of optimization is to minimize the compliance of the entire structure within the corresponding design domain.  The single material's Young's modulus is specified as $E_1 = 1.0$ and Poisson ratio as $\nu = 0.3$. The design domain is always discretized into quadrilateral elements, but with different resolutions. In Eq. \eqref{eq:non-convex-to}, $f(\mathbf{u}, \bm{\rho}) = \mathbf{f}^\intercal \mathbf{u}$, where $\mathbf{f}$ is the discretized external force vector; and there is only one constraint, i.e., the volume constraint that constrains the volume fraction of solid material toward the target value $V_T$. In each problem, we consider three different target volume fractions, i.e., $V_T=$ 0.3, 0.4, or 0.5, for a systematic analysis and validation. 
Particular emphasis is placed on the scenario with the highest discretization resolution and lowest target volume fraction, because of the high fidelity of designs provided by high discretization resolutions and the preference of low volume fractions of solid material in practical designs \cite{aage2017giga}. The solution process consistently started with an initial gray configuration satisfying the desired target volume fraction, same as that in SIMP \cite{andreassen2011efficient}. The convergence tolerance (as delineated in \Cref{algo:multi-cuts}) was set as $\varepsilon = 5 \times 10^{-3}$.

\subsubsection{Fixed trust-region radius} \label{subsec:single_minimum_compliance_fixedradius}
A fixed trust-region radius was firstly employed without dynamically varying it in the solution procedure. Thus, the trust-region radius was always equal to its initial value $d^0$, and three different fixed radii were examined, each with $d^0=0.05, 0.1$ or 0.2, respectively. In the parameter relaxation scheme proposed in \S \ref{subsec:parameter_relaxation}, only the minimum Young's modulus $E_0$ was relaxed. Thereby, the optimization procedure was split into two stages: in the first stage, it was set as $E_0 = 10^{-2}$; in the second stage, it was set with the required value in the original problem as $E_0 = 10^{-9}$. The numerical results are summarized in Table \ref{tab:fixed_radius_minimum_compliance}, where $N_{\text{FEM}}$ denotes the number of FEM analyses performed, or equivalently the number of iterations required throughout the optimization procedure; $f$ denotes the resultant value of the objective function. Here, the maximum allowed number of iterations was set as 100.
\begin{table}[ht]
    \caption{Single-material minimum compliance: Results with fixed trust-region radii.}
    \label{tab:fixed_radius_minimum_compliance}
    \centering
    \begin{tabular}{ccccccccccc}
         \toprule
         \multirow{2}{*}{$r$} & \multicolumn{5}{c}{\textbf{MBB}} & \multicolumn{5}{c}{\textbf{Cantilever}} \\
         \cmidrule(lr){2-6} \cmidrule(lr){7-11}
         & Resolution & $V_T$ & $d^0$ & $N_{\text{FEM}}$ & $f$ & Resolution & $V_T$ & $d^0$ & $N_{\text{FEM}}$ & $f$ \\
         \midrule
         \multirow{9}{*}{4} & \multirow{9}{*}{$240 \times 80$} & \multirow{3}{*}{0.3} & 0.05 & 100 & 446.36 & \multirow{9}{*}{$240 \times 120$} & \multirow{3}{*}{0.3} & 0.05 & 24 & 103.14 \\
         & & & 0.1 & 100 & 368.66 & & & 0.1 & 31 & 135.55 \\
         & & & 0.2 & 100 & 446.36 & & & 0.2 & 100 & 349.87 \\
         \cmidrule(lr){3-6} \cmidrule(lr){8-11}
         & & \multirow{3}{*}{0.4} & 0.05 & 34 & 241.60 & & \multirow{3}{*}{0.4} & 0.05 & 23 & 82.86 \\
         & & & 0.1 & 37 & 232.20 & & & 0.1 & 28 & 86.81 \\
         & & & 0.2 & 66 & 275.53 & & & 0.2 & 29 & 75.51 \\
         \cmidrule(lr){3-6} \cmidrule(lr){8-11}
         & & \multirow{3}{*}{0.5} & 0.05 & 30 & 210.79 & & \multirow{3}{*}{0.5} & 0.05 & 26 & 62.92 \\
         & & & 0.1 & 32 & 196.96 & & & 0.1 & 24 & 63.27 \\
         & & & 0.2 & 100 & 265.65 & & & 0.2 & 25 & 64.35 \\
         \midrule
         \multirow{9}{*}{6} & \multirow{9}{*}{$360 \times 120$} & \multirow{3}{*}{0.3} & 0.05 & 65 & 380.72 & \multirow{9}{*}{$360 \times 180$} & \multirow{3}{*}{0.3} & 0.05 & 100 & 115.12 \\
         & & & 0.1 & 100 & 442.64 & & & 0.1 & 33 & 121.38 \\
         & & & 0.2 & 100 & 742.72 & & & 0.2 & 100 & 381.12 \\
         \cmidrule(lr){3-6} \cmidrule(lr){8-11}
         & & \multirow{3}{*}{0.4} & 0.05 & 34 & 240.08 & & \multirow{3}{*}{0.4} & 0.05 & 24 & 82.93 \\
         & & & 0.1 & 100 & 251.10 & & & 0.1 & 27 & 77.85 \\
         & & & 0.2 & 100 & 355.44 & & & 0.2 & 35 & 79.22 \\
         \cmidrule(lr){3-6} \cmidrule(lr){8-11}
         & & \multirow{3}{*}{0.5} & 0.05 & 31 & 227.37 & & \multirow{3}{*}{0.5} & 0.05 & 26 & 62.92 \\
         & & & 0.1 & 30 & 196.03 & & & 0.1 & 24 & 63.59 \\
         & & & 0.2 & 29 & 192.64 & & & 0.2 & 25 & 64.97 \\
         \bottomrule
    \end{tabular}
\end{table}

From the results, it is obvious that using a fixed trust-region radius may not ensure reaching a converged solution within the maximum allowed number of iterations, especially when the target volume fraction is small (e.g. $V_T = 0.3$) or the trust-region radius is set too large (e.g. with $d^0 = 0.2$). This implies that if a fixed trust-region radius is employed, how to select its proper value should be problem dependent, and finding a consistent value applicable across different problem setups, each characterized by different BCs, discretization resolutions, and/or volume constraints, can be inherently challenging. Generally speaking, the optimization process with either too small or too large trust-region radii could require more iteration steps to reach convergence \cite{nocedal1999numerical}. More sophisticate discussions are provided below. 

As indicated in sensitivity analysis about the discrete TO problems \cite{sun2022sensitivity}, the accuracy of the sensitivities estimated for void elements ($\rho_i = 0$) is usually lower than that for solid elements ($\rho_i = 1$). Therefore, if the target volume fraction is large (e.g. $V_T = 0.5$), the fidelity of the sensitivity analysis can be maintained, and the proposed framework can easily converge within a reasonable number of iterations (around 30) with the optimal values of the objective function comparable across various choices of $d^0$.  On the contrary, if the target volume fraction is small (e.g. $V_T = 0.3$), a smaller trust-region radius (e.g., with $d^0 = 0.05$) is preferred, in order to enhance the sensitivity analysis' accuracy and ensure a fast converged solution. A representative instance matching this discussion is the MBB problem solved at the resolution of $360 \times 120$ for $V_T = 0.3$ or 0.4. The convergence of solution is achieved solely when $d^0$ is set to 0.05, resulting in 65 or 34 iteration counts, whereas setting $d^0$ to 0.1 or 0.2 both fails to yield convergence within the maximum allowed number of iterations. 
This suggests that if a fixed trust-region radius is employed, its value should be set sufficiently small to ensure convergence. However, if its value is set too small, the optimization process may be terminated early by the stopping criterion $\tfrac{|\eta^k - U|}{|U|} < \varepsilon$ before the local minimum of the objective function is actually found. For instance, 
in the MBB problem at the resolution of $360 \times 120$ and with the target volume fraction of $V_T = 0.5$, employing a small trust-region radius ($d^0 = 0.05$) results in a higher objective function value compared to using $d^0 = 0.2$, with the difference between the resultant objective function values as much as 18.0\%.

The findings of this study indicate that trust-region radii need to be determined by the problem and dynamically adjusted during the optimization process to ensure fast convergence of the solution.

\subsubsection{Adaptive trust-region radius}
\label{subsec:num-ex-adaptive-trust-region}
\begin{figure*}[ht]
    \centering
    \includegraphics[width=.8\textwidth]{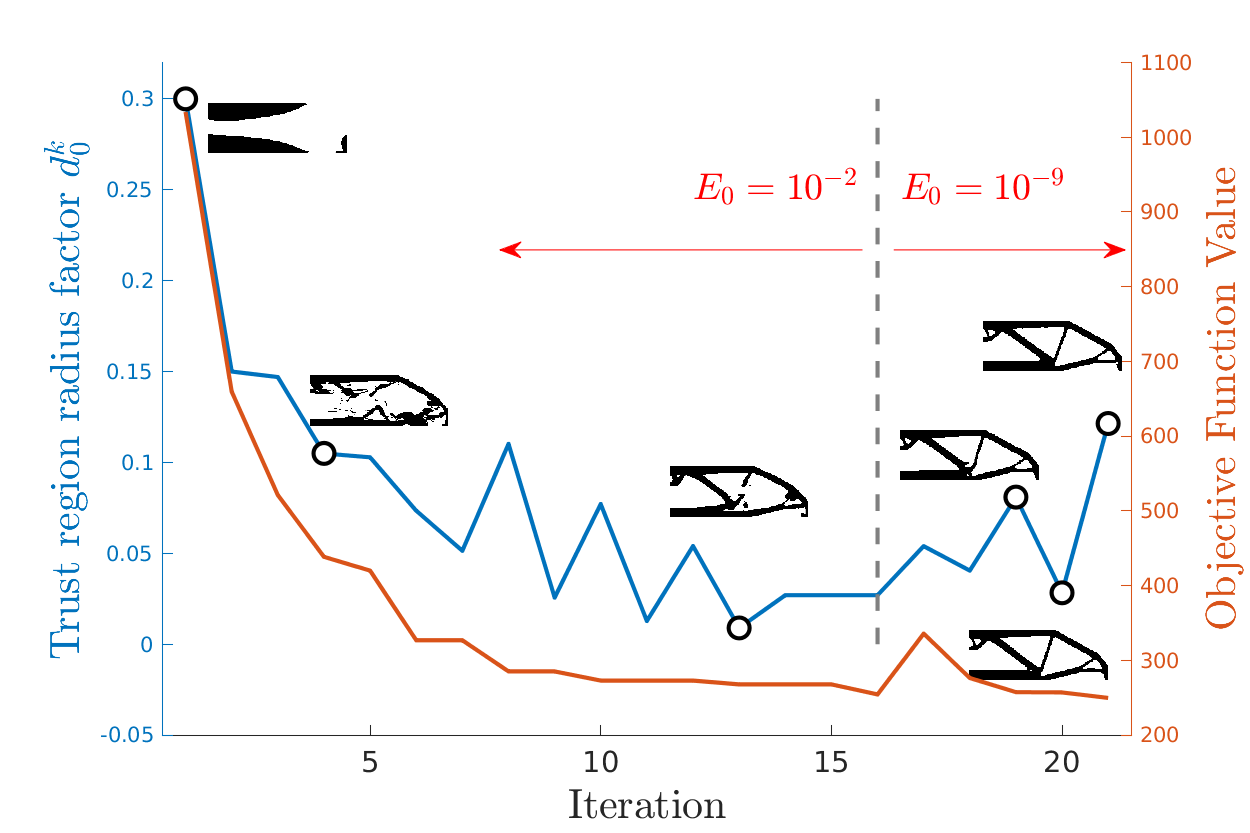}
    \caption{Single-material minimum compliance for the MBB design: Evolution of the adaptive trust-region radius and the objective function value during the optimization process, along with the resulting material configuration at different iteration steps. Here, the target volume fraction is $V_T = 0.4$; the discretization resolution is $240 \times 80$; and the initial trust-region radius is set as $d^0=0.3$. The minimum Young's moduli are $E_0 = 10^{-2}$ and $E_0 = 10^{-9}$, respectively, for the iteration steps before and after the gray dashed line.}
    \label{fig:iteration_TR_radius}
\end{figure*}

For the tests in this section, instead of maintaining a fixed value, the trust-region radius was dynamically varied following the rule described in Eq. \eqref{eq:trust_region_factor}. \Cref{fig:iteration_TR_radius} depicts the evolution of the trust-region radius and the objective function value during the optimization process, along with the resulting material configuration at different iteration steps for the MBB design.
Same as in the prior section, only the minimum Young's modulus $E_0$ was relaxed, and the optimization process was split into two stages with $E_0 = 10^{-2}$ in the first stage and $E_0 = 10^{-9}$ in the second stage, as indicated in \Cref{fig:iteration_TR_radius}.

The results for all the cases considered are compared to the best results achieved using fixed trust-region radii, as presented in \Cref{tab:fixed_adaptive_radius}. The column of \textbf{Fixed Radius} denotes the best performed cases with respect to the resultant objective function value seen in \Cref{tab:fixed_radius_minimum_compliance}. For the results using adaptive radii, as listed in the column of \textbf{Adaptive Radius}, $d^0$ denotes the \textit{initial} trust-region radius, whose value is then varied according to Eq. \eqref{eq:trust_region_factor} in subsequent iterations throughout the solution procedure. 
\begingroup
\setlength{\tabcolsep}{8pt} 
\begin{table*}[ht]
    \caption{Single-material minimum compliance: Comparison for employing adaptive vs. fixed trust-region radii.}
    \label{tab:fixed_adaptive_radius}
    \centering
    \resizebox{\textwidth}{!}{
    \begin{tabular}{ccccccc|ccc}
        \toprule
        \multirow{2}{*}{\textbf{Problem}} & \multirow{2}{*}{\textbf{Resolution}} & \multirow{2}{*}{$r$} & \multirow{2}{*}{$V_T$} & \multicolumn{3}{c}{\textbf{Fixed Radius}} & \multicolumn{3}{c}{\textbf{Adaptive Radius}} \\
        \cmidrule(lr){5-7} \cmidrule(lr){8-10} 
        & & & & $d^0$ & $N_{\text{FEM}}$ & $f$ & $d^0$ & $N_{\text{FEM}}$ & $f$ \\
        \midrule
        \multirow{6}{*}{MBB} & \multirow{3}{*}{$240 \times 80$} & \multirow{3}{*}{$4$} & $0.3$ & 0.1 & 100 & 368.66 & 0.3 & 36 & 294.43 \\
        & & & $0.4$ & 0.1 & 37 & 232.20 & 0.4 & 19 & 233.80 \\
        & & & $0.5$ & 0.1 & 32 & 196.96 & 0.4 & 20 & 193.45 \\
        \cmidrule(lr){2-10}
        & \multirow{3}{*}{$360 \times 120$} & \multirow{3}{*}{$6$} & $0.3$ & 0.05 & 65 & 380.72 & 0.4 & 28 & 315.02 \\
        & & & $0.4$ & 0.05 & 34 & 240.08 & 0.4 & 22 & 236.72 \\
        & & & $0.5$ & 0.2 & 29 & 192.64 & 0.4 & 26 & 214.60 \\
        \midrule
        \multirow{6}{*}{Cantilever} & \multirow{3}{*}{$240 \times 120$} & \multirow{3}{*}{$4$} & $0.3$ & 0.05 & 24 & 103.14 & 0.5 & 19 & 99.62 \\
        & & & $0.4$ & 0.2 & 29 & 75.51 & 0.5 & 15 & 77.79 \\
        & & & $0.5$ & 0.05 & 26 & 62.92 & 0.5 &  17 & 63.62 \\
        \cmidrule(lr){2-10}
        & \multirow{3}{*}{$360 \times 180$} & \multirow{3}{*}{$6$} & $0.3$ & 0.1 & 100 & 115.12 & 0.5 & 19 & 104.71 \\
        & & & $0.4$ & 0.1 & 27 & 77.85 & 0.5 & 20 & 76.54 \\
        & & & $0.5$ & 0.05 & 26 & 62.92 & 0.5 & 16 & 64.46 \\
        \hline
    \end{tabular}
    }
\end{table*}
\endgroup
For both MBB and cantilever problems and across various discretization resolutions and target volume fractions, employing adaptive trust-region radii consistently
ensured convergence to the optimal solution within the maximum allowed number of iterations. Next, the number of optimization iterations (or $N_{\text{FEM}}$) required for convergence remains notably small and relatively consistent across different scenarios. When compared to the outcomes obtained using fixed trust-region radii, $N_{\text{FEM}}$ is consistently smaller in each scenario considered. Finally, with respect to the resultant objective function value, while employing fixed radii may occasionally lead to slightly lower objective function values, particularly for higher target volume fractions (e.g., $V_T = 0.5$), the superiority of adaptive radii becomes apparent for lower target volume fractions (especially, $V_T = 0.3$), which are more relevant in practical designs. This superiority is characterized by not only achieving lower objective function values but also requiring fewer iteration counts, thus highlighting the effectiveness of adaptive radii over fixed radii when dealing with designs for lower target volume fractions.

Regarding the initial value of the trust-region radius $d^0$, we tested it ranging from 0.05 to 0.5. The best results with respect to $N_{\text{FEM}}$ are included in \Cref{tab:fixed_adaptive_radius}. It is worth noting that by adaptively adjusting the trust-region radius instead of using a fixed value, we could begin with a larger $d^0$ while achieving convergence and comparable objective function values with fewer iterations. And this becomes particularly evident when the target volume fraction is small (e.g. $V_T = 0.3$). In such cases, employing adaptive trust-region radii consistently yield lower objective function values and required far fewer iteration steps compared to employing fixed trust-region radii.

\subsubsection{Comparison with other TO methods}

In this section, we focus on the comparison with other TO methods, including SIMP with Heaviside projection \cite{sigmund2007morphology} and GBD \cite{munoz2011generalized,ye2023quantum}, for which the MBB problem was considered. The comparison results are summarized in \Cref{tab:comparison_minimum_compliance}.
\begin{table*}[ht]
    \caption{Single-material minimum compliance: Comparison with other TO methods for solving the MBB problem.}
    \label{tab:comparison_minimum_compliance}
    \centering
    \resizebox{\textwidth}{!}{
    \begin{tabular}{ccccc|cc|cc|cccc}
        \toprule
        \multirow{3}{*}{\textbf{Resolution}} & \multirow{3}{*}{$r$} & \multirow{3}{*}{$V_T$} & \multicolumn{2}{c}{\textbf{Our Method}} & \multicolumn{2}{c}{\textbf{SIMP-Heaviside}} & \multicolumn{2}{c}{\textbf{GBD-$E_0$}} & \multicolumn{4}{c}{\textbf{GBD-$V_T$}} \\
        \cmidrule(lr){4-5} \cmidrule(lr){6-7} \cmidrule(lr){8-9} \cmidrule(lr){10-13}
        & & & \multirow{2}{*}{$N_{\text{FEM}}$} & \multirow{2}{*}{$f$} & \multirow{2}{*}{$N_{\text{FEM}}$} & \multirow{2}{*}{$f$} & \multirow{2}{*}{$N_{\text{FEM}}$} & \multirow{2}{*}{$f$} & \multicolumn{2}{c}{$\Delta V = \frac{1}{24}$} & \multicolumn{2}{c}{$\Delta V = \frac{1}{12}$} \\
        \cmidrule(lr){10-11} \cmidrule(lr){12-13}
        & & & & & & & & & $N_{\text{FEM}}$ & $f$ & $N_{\text{FEM}}$ & $f$ \\
        \midrule
        \multirow{3}{*}{$240 \times 80$} & \multirow{3}{*}{$4$} & $0.3$ & 36 & 294.43 & 300 & 310.24 & 300 & 8.19e+10 & 119 & 286.95 & 104 & 414.51 \\
        & & $0.4$ & 19 & 233.80 & 300 & 239.53 & 300 & 5.80e+10 & 107 & 223.50 & 99 & 266.23 \\
        & & $0.5$ & 20 & 193.45 & 300 & 193.10 & 300 & 2.05e+10 & 82 & 193.18 & 76 & 190.07 \\
        \midrule
        \multirow{3}{*}{$360 \times 120$} & \multirow{3}{*}{$6$} & $0.3$ & 28 & 350.48 & 300 & 312.31 & 300 & 7.38e+10 & 104 & 454.84 & 72 & 496.75 \\
        & & $0.4$ & 22 & 236.72 & 300 & 238.30 & 300 & 6.18e+10 & 92 & 233.25 & 70 & 232.64 \\
        & & $0.5$ & 26 & 214.60 & 300 & 194.81 & 300 & 2.02e+10 & 76 & 190.92 & 54 & 194.25 \\
        \bottomrule
    \end{tabular}
    }
\end{table*}
Here, the column \textbf{Our Method} contains the results presented in the rightmost column of \Cref{tab:fixed_adaptive_radius}, i.e., with adaptive trust-region radii, for the MBB problem. For GBD, we examined two ways of parameter relaxation: one for the minimum Young's modulus and the other for the target volume fraction. In the former, the optimization procedure was split into two stages: $E_0 = 10^{-2}$ was employed in the first stage, and in the second stage, it was set with the required value in the original problem as $E_0 = 10^{-9}$, same as that employed in the prior two sections. The corresponding results are presented under the column \textbf{GBD-$E_0$} in \Cref{tab:comparison_minimum_compliance}. In the latter, the target volume fraction in the volume constraint is decreased incrementally per stage until reaching the desired value in the original problem, with the initial value set to 1.0 and the change increment fixed at $\Delta V = \frac{1}{24}$ or $\Delta V = \frac{1}{12}$. The corresponding results are presented under the column \textbf{GBD-$V_T$} in \Cref{tab:comparison_minimum_compliance}. The maximum number of iterations was capped at 300 for all the methods compared herein.   

Across different discretization resolutions and target volume fractions, our method consistently converges faster with significantly fewer iteration steps—about one order of magnitude fewer compared to the iteration counts needed by SIMP, as indicated by $N_{\text{FEM}}$, the number of FEM analyses called during the solution process. The resultant values of the objective function are comparable to, and sometimes even lower than, those obtained by the SIMP method, particularly for cases with lower target volume fractions. Our method also outperforms the GBD method in both the iteration counts and solution quality. When the minimum Young's modulus was chosen as the parameter to relax, GBD struggles to identify the optimal solutions, particularly for the cases with lower target volume fractions. Considering that the minimum compliance problem is a convex problem, solving it using the GBD method is equivalent to using our framework with a fixed trust-region radius of $d^k \equiv 1$. As discussed in \S\ref{subsec:single_minimum_compliance_fixedradius}, a too large fixed trust-region radius cannot ensure a stable and converged solution. When the target volume fraction (in the volume constraint) was chosen as the parameter to relax, GBD performed better. This is because by gradually changing the volume fraction, a varying trust region is implicitly enforced. The change in the design variable $\bm{\rho}$ in each iteration is constrained by the volume fraction increment $\Delta V$, which limits the allowable total change in $\bm{\rho}$. To ensure convergence, a smaller $\Delta V$ is preferred, but it could result in more iteration steps, which can be found when comparing the results for $\Delta V = \frac1{24}$ vs. $\Delta V = \frac1{12}$. 

\begin{figure*}[ht]
    \centering
    \begin{subfigure}[t]{.3\textwidth}
        \centering
        \includegraphics[width=.9\textwidth]{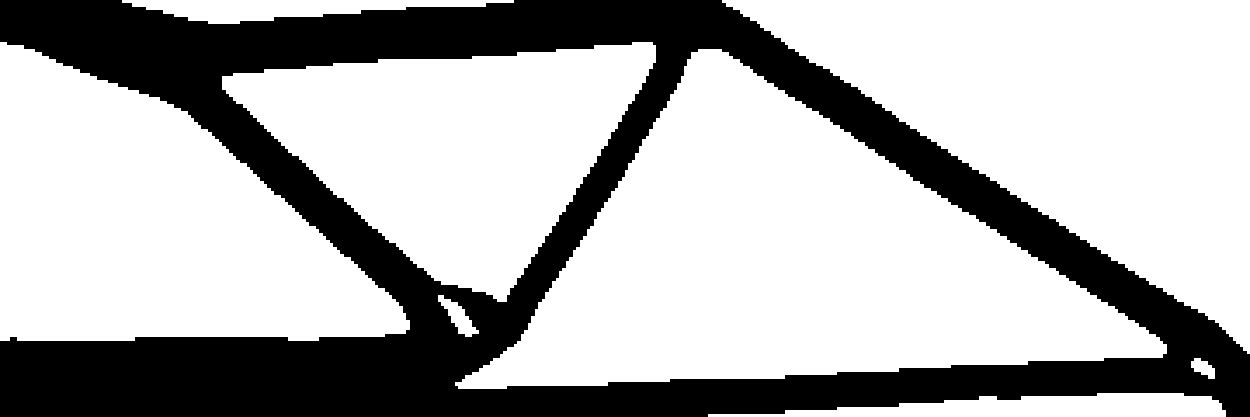}
        \caption{Our method}
        \label{subfig:multi-cuts-mc}
    \end{subfigure}
    \quad
    \begin{subfigure}[t]{.3\textwidth}
        \centering
        \includegraphics[width=.9\textwidth]{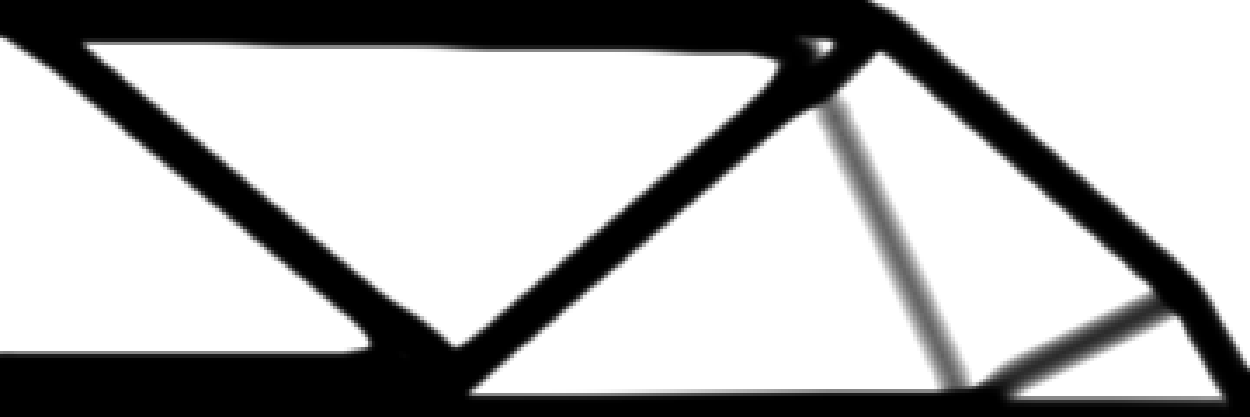}
        \caption{SIMP-Heaviside}
        \label{subfig:simp-mc}
    \end{subfigure}
    \quad
    \begin{subfigure}[t]{.3\textwidth}
        \centering
        \includegraphics[width=.9\textwidth]{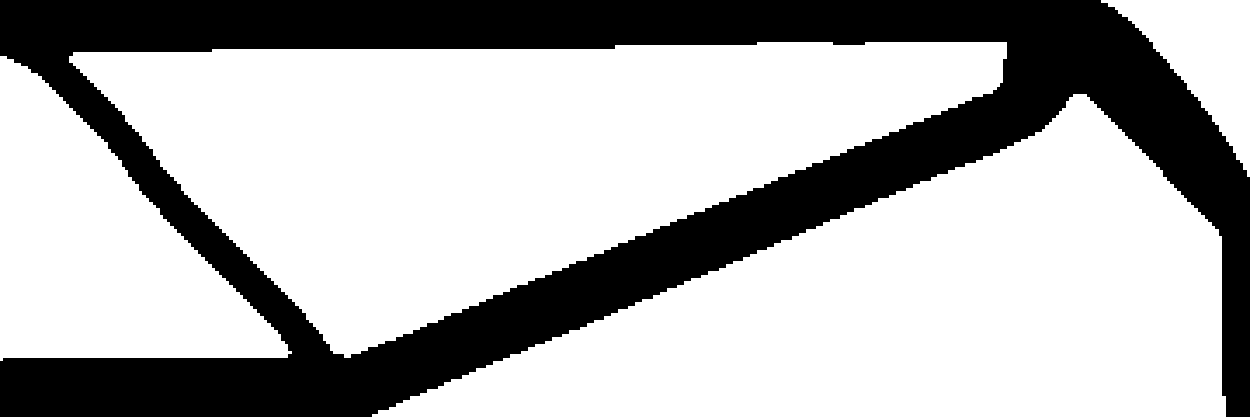}
        \caption{GBD-$V_T$}
        \label{subfig:gbd_minimum_compliance}
    \end{subfigure}
    \caption{Single-material minimum compliance: The optimal topology obtained for the MBB design from different methods, with the discretization resolution of $360 \times 120$ and the target volume fraction of $V_T = 0.3$.}
    \label{fig:minimum_compliance_360x120_0.3}
\end{figure*}

The resulting optimal topology is illustrated in \Cref{fig:minimum_compliance_360x120_0.3}. Our method achieved a clear-cut 0/1 material layout, in contrast to that produced by the SIMP method, where part of the topology is still in gray scale, requiring more optimization iterations or post-processing steps to reach a clear-cut 0/1 configuration. For \textbf{GBD-$V_T$} with $\Delta V = \frac1{24}$, although the optimization converged, the resulting objective function value remained high, higher than those obtained by our method and SIMP by more than 30\%. This discrepancy leads to a notably distinct final topology generated by the GBD method, as depicted in \Cref{fig:minimum_compliance_360x120_0.3}.

However, the resulting optimal topology in \Cref{subfig:multi-cuts-mc} with the deflected top-left and bottom-right parts, may be still undesirable. To further improve the topology, the parameter relaxation scheme (proposed in \S\ref{subsec:parameter_relaxation}) was employed, as discussed in the next subsection. The  improved topology is then shown in \Cref{subfig:scheme2_360x120}.

\subsubsection{Parameter relaxation}
\label{subsubsec:parameter-adjustment-improvement}
In this section, we present a thorough analysis for parameter relaxation, exploring its effectiveness in improving both solution efficiency and quality. We discuss the rationale behind its efficacy and the selection of suitable parameters for relaxation. 

During the solution process, especially in the initial stages, the optimization problem encountered may suffer from ill-conditioning. This implies that slight variations on the design variable (i.e., a few elements transitioning from solid to void or vice versa), could lead to significant fluctuations in the objective function value. To demonstrate this, we first compare the objective function values corresponding to different topological configurations obtained with 
the minimum Young's modulus $E_0 = 10^{-9}$ vs. $E_0 = 10^{-2}$, as denoted as $f_1$ and $f_2$, respectively, in Table. \ref{tab:min_youngs}. 
\begin{table*}[ht]
    \centering
    \caption{Single-material minimum compliance for the MBB design: The roles of the minimum Young's modulus and the target volume fraction in the optimization problem's conditioning.}
    \label{tab:min_youngs}
    \begin{tabular}{N|c|c|c}
        \toprule
        \multicolumn{1}{c|}{Configuration} & $V_T$ & $f_1~ (E_0 = 10^{-9})$ & $f_2~(E_0 = 10^{-2})$ \\
        \midrule
        \label{fig:gbd2}
        \begin{minipage}{.3\textwidth}
            \includegraphics[width=\textwidth]{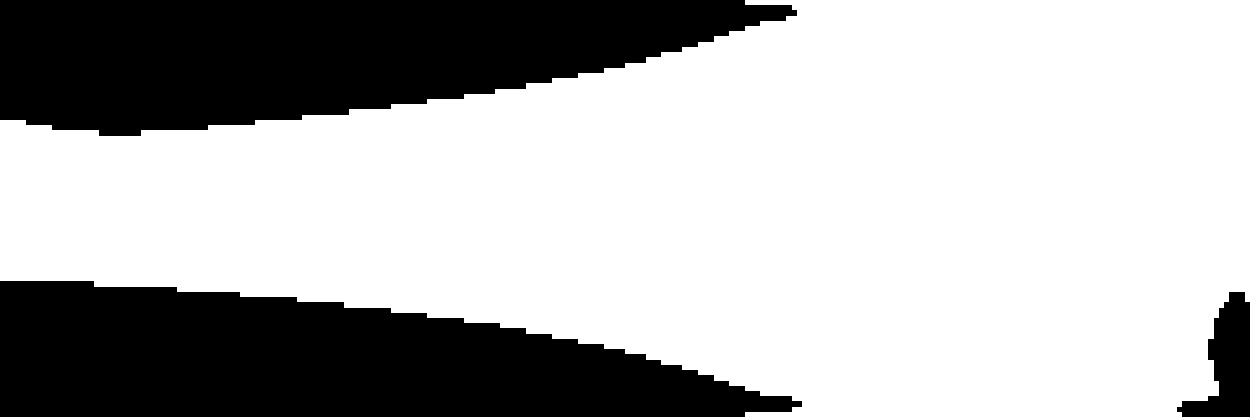}
        \end{minipage} & 0.3 & $10186564604.77$ & $1473.14$ \\
        \midrule
        \label{fig:gbd3}
        \begin{minipage}{.3\textwidth}
            \includegraphics[width=\textwidth]{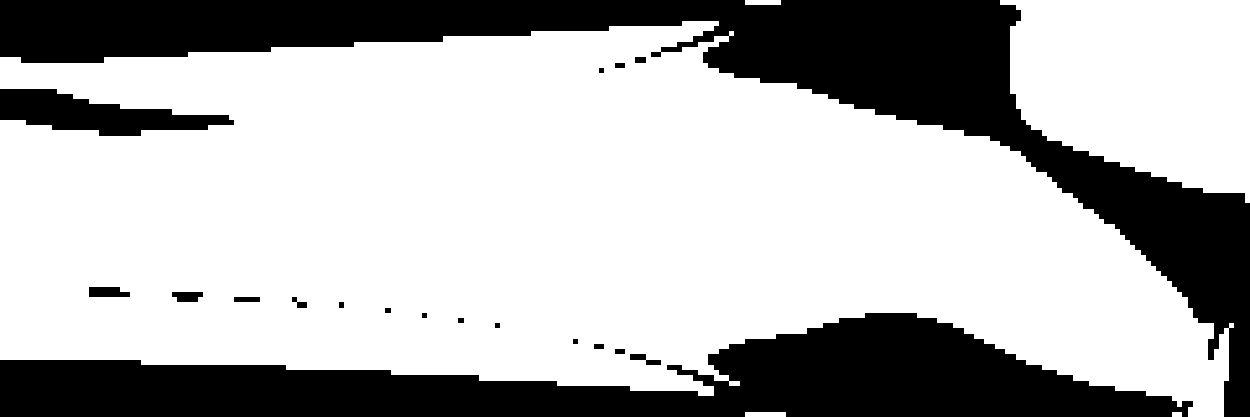}
        \end{minipage} & 0.3 & $109405.30$ & $764.62$ \\
        \midrule
        \label{fig:gbd2_0.6}
        \begin{minipage}{.3\textwidth}
            \includegraphics[width=\textwidth]{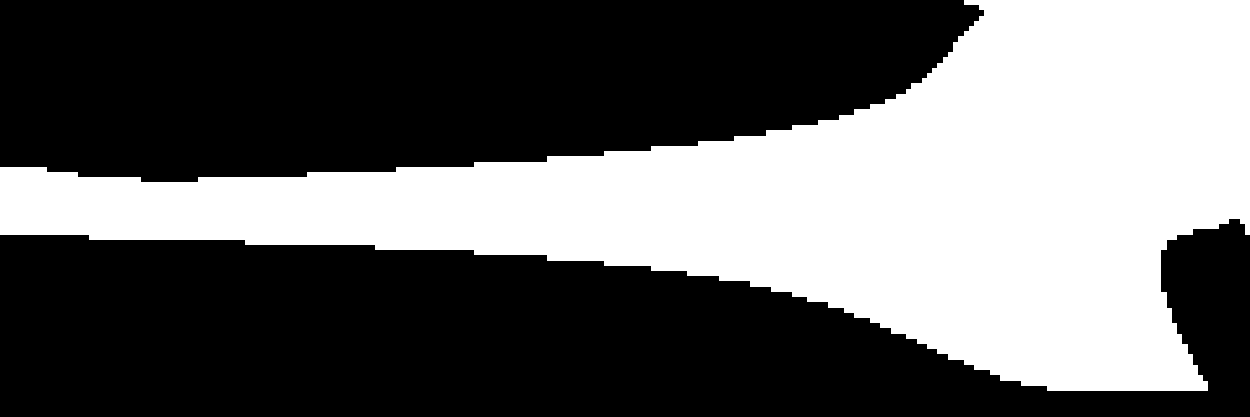}
        \end{minipage} & 0.6 & $83637287.04$ & $521.03$ \\
        \midrule
        \label{fig:gbd3_0.6}
        \begin{minipage}{.3\textwidth}
            \includegraphics[width=\textwidth]{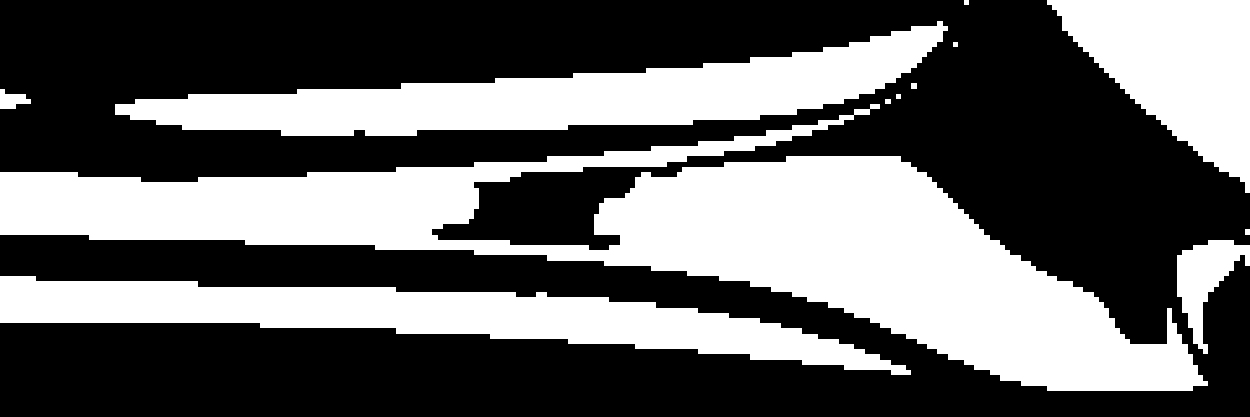}
        \end{minipage} & 0.6 & $2902.42$ & $398.42$ \\
        \bottomrule
    \end{tabular}
\end{table*}
Table \customref{tab:min_youngs}{fig:gbd2} and \customref{tab:min_youngs}{fig:gbd3} are what we obtained from the first two consecutive iteration steps, where the target volume fraction was set as $V_T = 0.3$. In these two steps, the difference in the design variable, measured in 2-norm, is $\frac{\| \bm{\rho}^1 - \bm{\rho}^2 \|_2}{n_e} = 0.3$, which is considered a minor change. However, the substantial disparity in objective function values between these two steps, with $E_0 = 10^{-9}$, indicates the presence of ill-conditioning in the optimization problem at this stage. Solving this ill-conditioned problem could lead to an unstable or non-converging solution. In contrast, by relaxing the minimum Young's modulus to $E_0 = 10^{-2}$, the ill-conditioning is mitigated. This elucidates why relaxing the value of the minimum Young's modulus during the optimization process can effectively enhance both solution efficiency and quality. We notice that in literature, other TO methods that rely on discrete programming chose to relax the target volume fraction during the optimization process, e.g., as reported in \cite{liang2019topology, picelli2021101, ye2023quantum}. Thus, we repeated the same numerical experiment but with a larger target volume fraction, i.e., $V_T = 0.6$, as outlined in Table \customref{tab:min_youngs}{fig:gbd2_0.6} and \customref{tab:min_youngs}{fig:gbd3_0.6}. The difference in the design variable $\bm{\rho}$, measured in 2-norm, is less than 0.3 between the two configurations in Table \customref{tab:min_youngs}{fig:gbd2_0.6} and \customref{tab:min_youngs}{fig:gbd3_0.6}. With $E_0 = 10^{-9}$, by setting a larger target volume fraction from 0.3 to 0.6, the optimization problem is still ill-conditioning. It is evident that enlarging the target volume fraction has a limited effect on improving the problem's conditioning, when compared with elevating the minimum Young's modulus. Thus, in the prior sections, we chose to focus on the minimum Young's modulus $E_0$ in parameter relaxation. 

\begin{figure*}[ht]
    \centering
    \begin{subfigure}[t]{.45\textwidth}
        \centering
        \includegraphics[height=1.5cm]{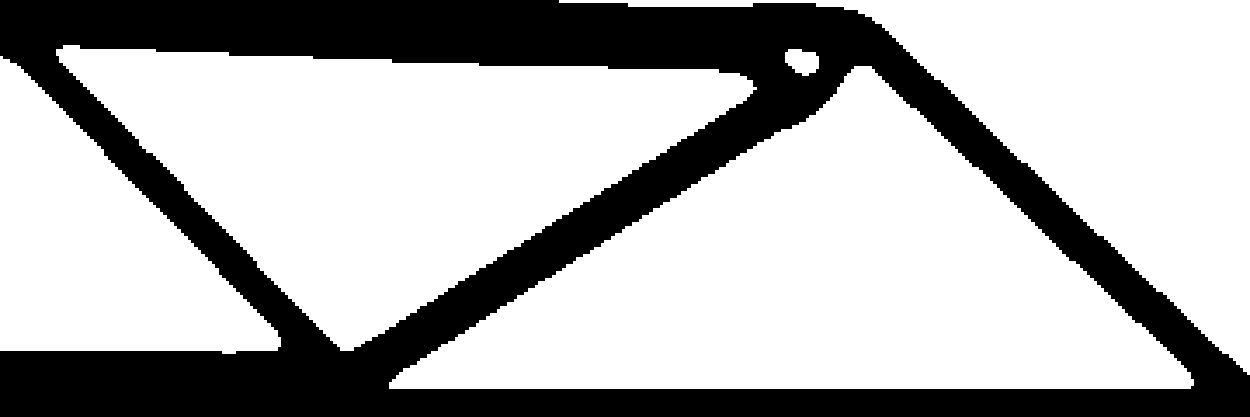}
        \caption{Length-height ratio: $3:1$; discretization resolution: $360 \times 120$, obtained by employing \textbf{Scheme 2} in \Cref{tab:MBB_relaxation_scheme} for parameter relaxation.}
        \label{subfig:scheme2_360x120}
    \end{subfigure}
    \quad \quad    
    \begin{subfigure}[t]{.45\textwidth}
        \centering
        \includegraphics[height=1.5cm]{Figures/top_360_120_03.png}
        \caption{Length-height ratio: $3:1$; discretization resolution: $360 \times 120$, obtained by SIMP}
    \end{subfigure}
    \begin{subfigure}[t]{.45\textwidth}
        \centering
        \includegraphics[height=1.5cm]{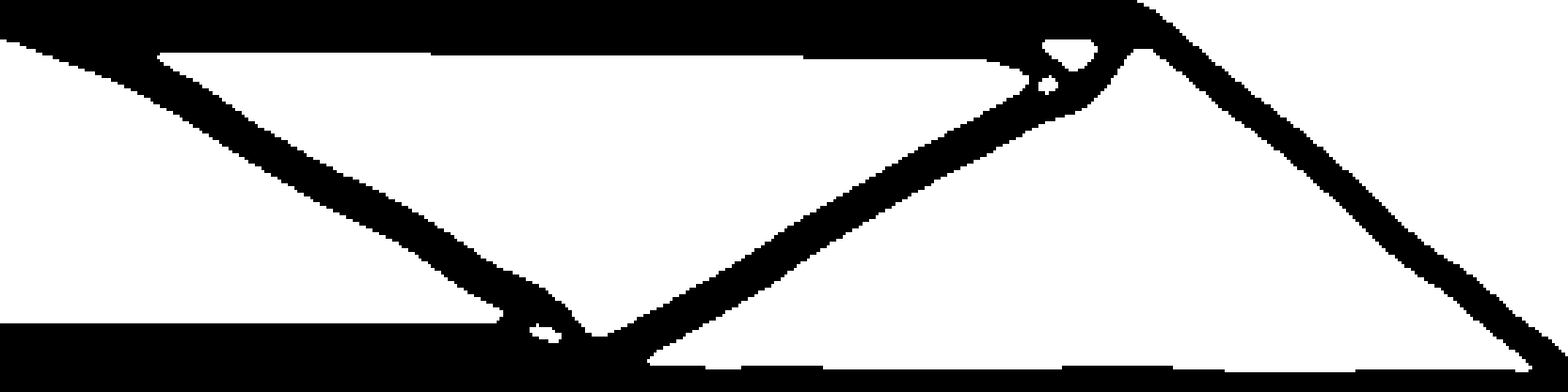}
        \caption{Length-height ratio: $4:1$; discretization resolution: $480 \times 120$, obtained by employing \textbf{Scheme 2} in \Cref{tab:MBB_relaxation_scheme} for parameter relaxation.}
    \end{subfigure}
    \quad \quad    
    \begin{subfigure}[t]{.45\textwidth}
        \centering
        \includegraphics[height=1.5cm]{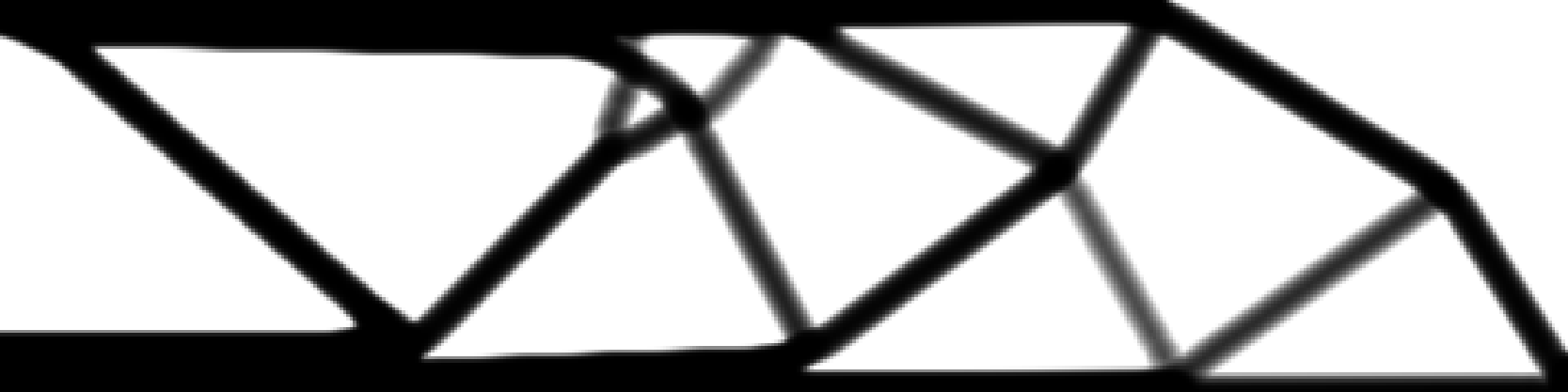}
        \caption{Length-height ratio: $4:1$; discretization resolution: $480 \times 120$, obtained by SIMP}
    \end{subfigure}
    \caption{Single-material minimum compliance for the MBB design: The final topological configuration with the target volume fraction of $V_T = 0.3$}
    \label{fig:parameter_relaxation_mbb}
\end{figure*}

In practice, relaxing more than one parameters can further improve the solution efficiency and/or solution quality, depending on the conditioning of the problem. Henceforth, we further explored adjusting both the minimum Young's modulus $E_0$ and the target volume fraction $V_T$. One exemplary implementation is given as \textbf{Scheme 2} in \Cref{tab:MBB_relaxation_scheme}, where the parameter relaxation process described in \S \ref{subsec:parameter_relaxation} was partitioned into 9 stages with $N_P=8$ and the initial value for the volume fraction $V_1$ set as 0.6 in Eq. \eqref{eq:relaxation_interpolation}. The results, as summarized in \Cref{tab:MBB_relaxation_scheme_result}, are compared with those obtained by relaxing only $E_0$ in two stages as employed in the prior sections, denoted as \textbf{Scheme 1} in \Cref{tab:MBB_relaxation_scheme}. 
\begin{table}[ht]
    \centering
    \caption{Single-material minimum compliance for the MBB design: Two parameter relaxation schemes.}
    \label{tab:MBB_relaxation_scheme}
    \begin{tabular}{ccccccccccc}
        \toprule
        & Stage & 1 & 2 & 3 & 4 & 5 & 6 & 7 & 8 & 9 \\
        \midrule
        \multirow{2}{*}{Scheme 1} & $V_T$ & 0.3 & 0.3 \\
        & $E_0$ & $10^{-2}$ & $10^{-9}$ \\
        \midrule
        \multirow{2}{*}{Scheme 2} & $V_T$ & 0.600 & 0.543 & 0.492 & 0.446 & 0.404 & 0.366 & 0.331 & 0.300 & 0.300 \\
        & $E_0$ & $10^{-2}$ & $10^{-2}$ & $10^{-2}$ & $10^{-2}$ & $10^{-2}$ & $10^{-2}$ & $10^{-2}$ & $10^{-2}$ & $10^{-9}$ \\
        \bottomrule
    \end{tabular}
\end{table}
For a more systematic analysis, we also extended the design domain to a larger length-height ratio, $L:H=4:1$. For the domain of each ratio, two different discretization resolutions were examined. The results are summarized in \Cref{tab:MBB_relaxation_scheme_result}. 
\begingroup
\setlength{\tabcolsep}{10pt} 
\begin{table*}[ht]
    \caption{Single-material minimum compliance for the MBB design: Comparison between the two parameter relaxation schemes, using the results from the SIMP method as the reference.}
    \label{tab:MBB_relaxation_scheme_result}
    \centering
    \begin{tabular}{ccc|cc|cc}
        \toprule
        \multirow{2}{*}{\textbf{Resolution}} & \multicolumn{2}{c}{\textbf{Scheme 1}} & \multicolumn{2}{c}{\textbf{Scheme 2}} & \multicolumn{2}{c}{\textbf{SIMP-Heaviside}} \\
        \cmidrule(lr){2-3} \cmidrule(lr){4-5} \cmidrule(lr){6-7}
        & $f$ & $N_{\text{FEM}}$ & $f$ & $N_{\text{FEM}}$ & $f$ & $N_{\text{FEM}}$ \\
        \midrule
        $240 \times 80$ & 294.43 & 36 & 294.07 & 52 & 310.24 & 300 \\
        $320 \times 80$ & 662.61 & 41 & 613.37 & 58 & 635.74 & 300 \\
        $360 \times 120$ & 350.48 & 28 & 292.88 & 51 & 312.31 & 300 \\
        $480 \times 120$ & 808.75 & 84 & 605.63 & 51 & 639.60 & 300 \\
        \bottomrule
    \end{tabular}
\end{table*}
\endgroup
The final material configuration for each domain, obtained with the finest discretization resolution, are exhibited in \Cref{fig:parameter_relaxation_mbb}. Through comparison, we find that the resultant objective function values from \textbf{Scheme 2} are consistently lower. More pronounced improvements are observed at higher discretization resolutions and greater aspect ratios of the design domain, e.g., in the case of $480\times 120$, with the improvement in the objective function value reaching as high as 26.4\% and the iteration counts reduced by about 40\%. These findings suggest that relaxing the minimum Young's modulus along with another parameter in multiple stages could potentially offer greater advantages in achieving lower objective function values and possibly reducing iteration counts as well. The extent of these enhancements can vary depending on the specific problem, particularly its conditioning.

\subsection{Single-material Compliant Mechanism}
\label{subsec:single-mat-cm}

We next move to a benchmark non-convex TO problem, the compliant mechanism design, where the displacement objective exhibits non-convexity with respect to the design variables $\bm{\rho}$. It represents a broader class of TO challenges than the minimum compliance design addressed earlier, which is regarded as a convex TO problem. 

\begin{figure*}[ht]
    \centering
    \includegraphics[width=.45\textwidth]{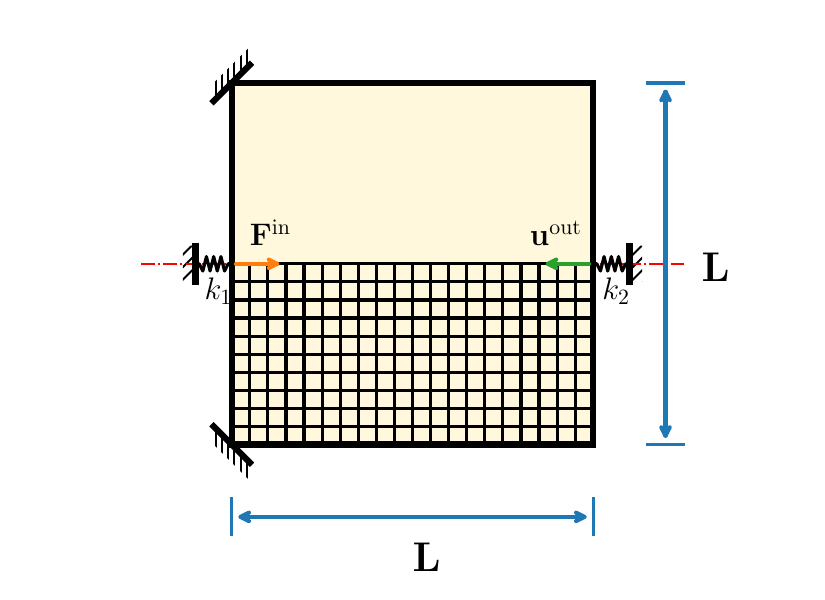}
    \caption{Single-material Compliant Mechanism: Design domain and boundary conditions.}
    \label{fig:compliant_mechanism}
\end{figure*}

The design domain and the associated BCs are depicted in \Cref{fig:compliant_mechanism}, where the domain is a square with each side measuring $L$, and the top-left and bottom-left corners are fixed in place. The domain is also connected to two springs each at the middle of the left or right side, respectively. The stiffness coefficient of both springs is $k_1 = k_2 = 0.1$. In addition, an input force $\mathbf{F}^{\text{in}}_x = 1$ is applied at the midpoint of the left side. The TO design herein aims to maximize the displacement along $x$-direction at the midpoint of the right side, i.e., minimizing $\mathbf{u}^{\text{out}}_x$ as denoted in \Cref{fig:compliant_mechanism}. Due to the problem's axial symmetry around the center line (as highlighted in red dotted-dash in \Cref{fig:compliant_mechanism}), only half of the design domain needs to be considered in the actual solution procedure and hence was discretized with a uniform mesh. The target volume fraction is $V_T = 0.3$. The material parameters are specified as: Young's modulus $E_1 = 1.0$; Poisson's ratio $\nu = 0.3$. The convergence tolerance required in \Cref{algo:multi-cuts} was set as $\varepsilon = 5 \times 10^{-3}$.

\begin{figure*}[ht]
    \centering
    \begin{subfigure}[t]{.45\textwidth}
        \centering
        \includegraphics[width=.6\textwidth]{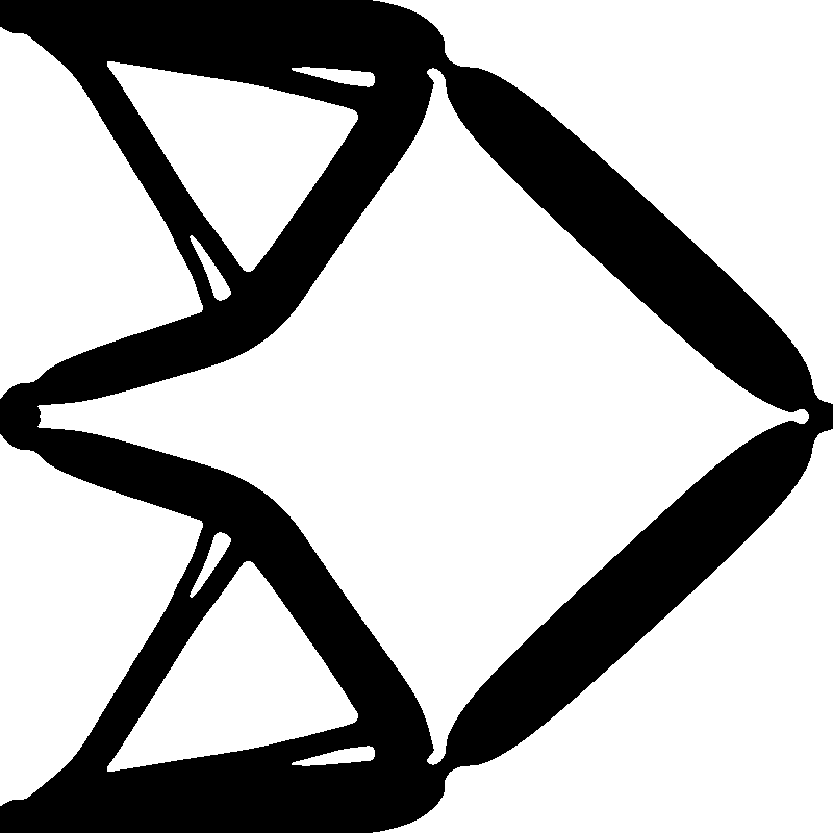}
        \caption{Our method with $d^0 = 0.3$}
        \label{fig:multicuts_cm}
    \end{subfigure}
    \quad \quad    
    \begin{subfigure}[t]{.45\textwidth}
        \centering
        \includegraphics[width=.6\textwidth]{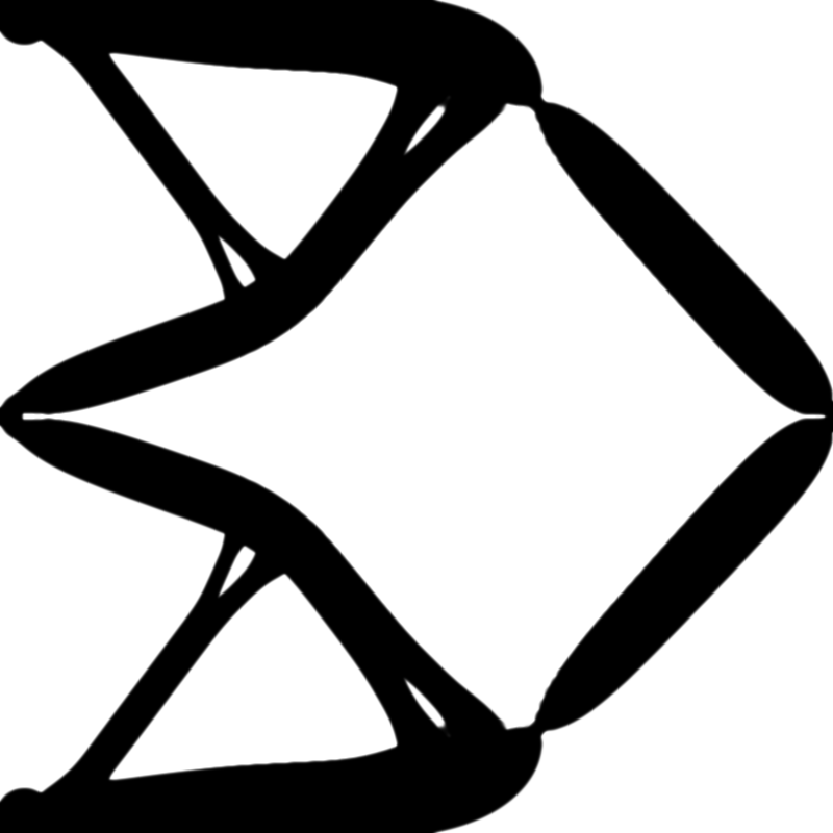}
        \caption{SIMP-Heavisde}
    \end{subfigure}
    \caption{Single-material compliant mechanism: The optimal topology obtained from different methods, with the discretization resolution of $800 \times 400$ and the target volume fraction of $V_T = 0.3$.}
    \label{fig:single_compliant_mechanism}
\end{figure*}

During the optimization process, the trust-region radii were adaptively adjusted following Eq. \eqref{eq:trust_region_factor}, same as that in \S \ref{subsec:num-ex-adaptive-trust-region}. Only the parameter $E_0$ was relaxed, and hence the optimization process was split into two stages: $E_0 = 10^{-2}$ was employed in the first stage; and it was set to the required value as $E_0 = 10^{-9}$ in the second stage. The results, including the iteration counts $N_{\text{FEM}}$ and the resultant objective function values $f$, are compared with those obtained from the SIMP method with Heaviside projection \cite{andreassen2011efficient}. The maximum number of iterations for the SIMP method was capped at 300. As indicated in \Cref{tab:comparison_compliant_mechanism}, across different discretization resolutions, our method, incorporating the multi-cut formulation along with adaptive trust-region radii and parameter relaxation, consistently achieved comparable or lower objective function values with significantly fewer iteration counts (about one order of magnitude fewer) compared to the SIMP method. Reducing $N_{\text{FEM}}$ by one order of magnitude can lead to significant savings in computational time in the execution of time-consuming FEM analyses, especially when tackling large-scale problems. The resultant optimal topology is depicted in \Cref{fig:single_compliant_mechanism}.
\begingroup
\setlength{\tabcolsep}{10pt} 
\begin{table*}[tp]
    \caption{Single-material compliant mechanism: Results and comparison with the SIMP method.}
    \label{tab:comparison_compliant_mechanism}
    \centering
    \begin{tabular}{ccccc|cc}
        \hline
        \multirow{2}{*}{\textbf{Resolution}} & \multirow{2}{*}{$r$} & \multicolumn{3}{c}{\textbf{Our Method}} & \multicolumn{2}{c}{\textbf{SIMP-Heaviside}} \\
        \cmidrule(lr){3-5} \cmidrule(lr){6-7}
        & & $d^0$ & $N_{\text{FEM}}$ & $f$ & $N_{\text{FEM}}$ & $f$ \\
        \midrule
        \multirow{6}{*}{$200 \times 100$} & \multirow{6}{*}{$2$} & 0.05 & 35 & -0.9271 & \multirow{6}{*}{$300$} & \multirow{6}{*}{-0.8960} \\
        & & 0.1 & 24 & -0.9139 & & \\
        & & 0.2 & 26 & -0.8919 & & \\
        & & 0.3 & 24 & -0.9226 & & \\
        & & 0.4 & 24 & -0.9098 & & \\
        & & 0.5 & 29 & -0.9282 & & \\
        \midrule
        \multirow{6}{*}{$400 \times 200$} & \multirow{6}{*}{$4$} & 0.05 & 29 & -0.8661 & \multirow{6}{*}{$300$} & \multirow{6}{*}{-0.8511} \\
        & & 0.1 & 25 & -0.8652 & & \\
        & & 0.2 & 36 & -0.8677 & & \\
        & & 0.3 & 31 & -0.8602 & & \\
        & & 0.4 & 38 & -0.8728 & & \\
        & & 0.5 & 33 & -0.8763 & & \\
        \midrule
        \multirow{6}{*}{$800 \times 400$} & \multirow{6}{*}{$8$} & 0.05 & 36 & -0.8313 & \multirow{6}{*}{$300$} & \multirow{6}{*}{-0.8132} \\
        & & 0.1 & 28 & -0.8238 & & \\
        & & 0.2 & 25 & -0.8056 & & \\
        & & 0.3 & 28 & -0.8167 & & \\
        & & 0.4 & 27 & -0.8085 & & \\
        & & 0.5 & 37 & -0.8270 & & \\
        \bottomrule
    \end{tabular}
\end{table*}
\endgroup

We also systematically examined the impact of the choice of initial trust-region radius $d^0$, spanning from 0.05 to 0.5. The comparison reveals that the objective function values attained  closely align across different initial trust-region radii, with the maximum relative deviation between the highest and the lowest $f$ being under 3.9\%. Additionally, the number of FEM analyses $N_{\text{FEM}}$ or iteration counts required to reach solution convergence tends to remain consistent across all scenarios, notably not exceeding 40. This consistency underscores the robustness of the proposed method and its insensitivity to the choice of initial trust-region radii. 

\subsection{Multi-material Minimum Compliance}
Following the completion of benchmark single-material TO problems, including both convex and non-convex cases, we proceeded to address more challenging TO problems involving multiple materials. When addressing the inclusion of multiple candidate materials in an TO framework, a vital consideration is whether the proposed method can accurately discern the coexistence of all or some of the candidate materials in the derived material configuration. We first tackled a problem involving two materials, from which we demonstrate the effectiveness of the proposed parameter relaxation scheme for further enhancing solution quality (in terms of the objective function value). Next, we considered a design scenario incorporating five different candidate materials, by which we aimed to evaluate our framework's ability to identify and select the optimal combination of candidate materials in the resulting material layout. We refer to two recent studies \cite{liu2024multi, huang2021new} to validate our results and assess the performance of our framework.

\subsubsection{Two-material MBB}
Our examination started with a MBB design involving two distinct materials. The problem setup was the same as depicted in \Cref{subfig:mbb}, with the length-height ratio of the design domain as $L:H=2:1$. The design domain was discretized by a uniform mesh with $120 \times 60$ quadrilateral elements. Accordingly, a filter radius of $r = 2$ was employed. Same as in the prior sections, the convergence tolerance was set as $\varepsilon = 5 \times 10^{-3}$. The two materials are differentiated by their Young's moduli (specified as $E_1 = 0.55$ and $E_2 = 1.0$, respectively) and the normalized densities (specified as $\bar{M}_1 = 0.5$ and $\bar{M}_2 = 1.0$, respectively). The maximum permissible total mass fraction was set to $\bar{M}_{\max} = 0.4$.

We examined two different parameter relaxation schemes, highlighting that refining the parameter relaxation can improve the solution quality in terms of the objective function's value. In particular, relaxing the maximum permissible total mass fraction $\bar{M}_{\max}$ in multiple stages in addition to $E_0$ can further improve the conditioning of the optimization problem, similarly to our finding about relaxing the target volume fraction as discussed in \S \ref{subsubsec:parameter-adjustment-improvement}. As outlined in \Cref{tab:two-material-relaxation}, \textbf{Scheme 1} (as a reference) involves only two stages: $E_0 = 10^{-2}$ is employed in the first stage, followed with $E_0 = 10^{-9}$ in the second stage; in both stages, $\bar{M}_{\max}$ is constantly set as the desired value $0.4$. In \textbf{Scheme 2}, the parameter relaxation is partitioned into 8 stages. In the first 7 stages, $\bar{M}_{\max}$ is progressively decreased from 0.6 (initial value) to 0.4 (desired value) following Eq. \eqref{eq:relaxation_interpolation} with $N_P=7$. In all these stages, the minimum Young's modulus $E_0$ is maintained at $10^{-2}$. In the last stage, $E_0$ is decreased to the desired value $E_0 = 10^{-9}$, while $\bar{M}_{\max}$ stays unchanged. 
\begin{table*}[ht]
    \caption{Two-material MBB: Two parameter relaxation schemes.}
    \label{tab:two-material-relaxation}
    \centering
    \begin{tabular}{cccccccccc}
        \toprule
        & Stage & 1 & 2 & 3 & 4 & 5 & 6 & 7 & 8 \\
        \midrule
        \multirow{2}{*}{Scheme 1} & $\bar{M}_{\max}$ & 0.40 & 0.40 \\
        & $E_0$ & $10^{-2}$ & $10^{-9}$ \\
        \midrule
        \multirow{2}{*}{Scheme 2} & $\bar{M}_{\max}$ & 0.600 & 0.561 & 0.524 & 0.490 & 0.458 & 0.428 & 0.400 & 0.400 \\
        & $E_0$ & $10^{-2}$ & $10^{-2}$ & $10^{-2}$ & $10^{-2}$ & $10^{-2}$ & $10^{-2}$ & $10^{-2}$ & $10^{-9}$ \\
        \bottomrule
    \end{tabular}
\end{table*}

The main results are compiled in \Cref{tab:two-material-mbb}, including the number of FEM analyses $N_{\text{FEM}}$ invoked before reaching convergence and the resultant objective function values $f$.
\begin{table*}[ht]
    \caption{Two-material MBB: Main results of our method, comparison between two parameter relaxation schemes and with the reference results reported in \cite{liu2024multi}.}
    \label{tab:two-material-mbb}
    \centering
    \begin{tabular}{cccccc}
        \toprule
        \multicolumn{2}{c}{\textbf{Our Method} (Scheme 1)} & \multicolumn{2}{c}{\textbf{Our Method} (Scheme 2)} & \multicolumn{2}{c}{\textbf{Reference} \cite{liu2024multi}} \\
        \cmidrule(lr){1-2} \cmidrule(lr){3-4} \cmidrule(lr){5-6}
        $N_{\text{FEM}}$ & $f$ & $N_{\text{FEM}}$ & $f$ & $N_{\text{FEM}}$ & $f$ \\
        \midrule
        13 & 93.3632 & 32 & 83.3298 & 224 & 83.43 \\
        \bottomrule
    \end{tabular}
\end{table*}
Our method achieved comparable resultant objective function values while requiring significantly fewer iteration counts (or the number of FEM analyses) compared to the results reported in the reference \cite{liu2024multi}. This efficiency translates into considerable computational savings without compromising the quality of the resultant material layout, as gauged by the objective function value. It hence demonstrates the superior efficacy of the proposed framework for multi-material TO design. 

\begin{figure*}[ht]
    \centering
    \begin{subfigure}[t]{.45\textwidth}
        \centering
        \includegraphics[width=0.75\textwidth]{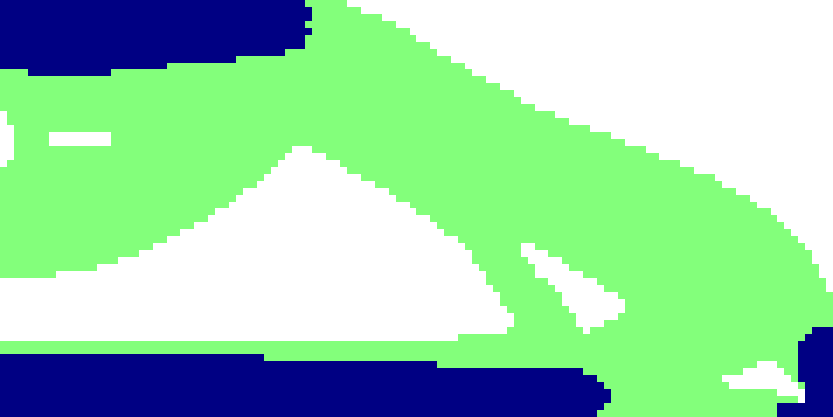}
        \caption{Scheme 1.}
    \end{subfigure}
    \quad \quad
    \begin{subfigure}[t]{.45\textwidth}
        \centering
        \includegraphics[width=0.75\textwidth]{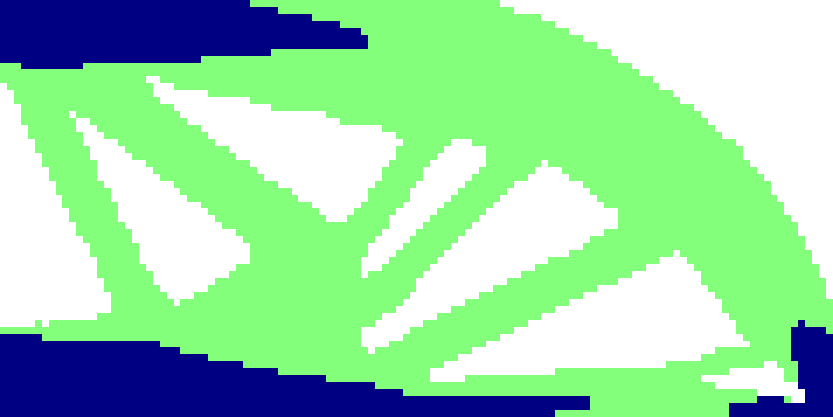}
        \caption{Scheme 2.}
    \end{subfigure}
    \caption{Two-material MBB: The optimal material configurations obtained by our method. Here, the discretization resolution is $120 \times 60$, and the maximum permissible total mass fraction is $\bar{M}_{\max} = 0.4$.  {\textcolor{mc_mat1}{\rule{2ex}{2ex}}} denotes MAT 1.  {\textcolor{mc_mat2}{\rule{2ex}{2ex}}} denotes MAT 2.}
    \label{fig:two-material-mbb}
\end{figure*}

As advocated in preceding sections, maintaining a relatively large value for $E_0$ until the concluding stage promotes the stability of the optimization procedure. Further refining the parameter relaxation, particularly for the maximum permissible total mass fraction $\bar{M}_{\max}$, can facilitate improving the solution quality for multi-material TO design, as indicated by the resultant objective function value. The material configurations derived from the two distinct parameter relaxation schemes are presented in \Cref{fig:two-material-mbb}. Despite the application of distinct parameter relaxation schemes, the core structures of the resultant material layouts remain consistent. And the distribution and placement of the two candidate materials also demonstrate considerable agreement across the results.
This indicates that employing parameter relaxation in different ways does not compromise the integrity of the primary structural topology optimized by the proposed TO framework; it merely refines the detailed structure to optimize the objective function value. 

\subsubsection{Five-material Cantilever}
We next tackled a cantilever design problem involving five materials to further assess the applicability and effectiveness of our method for general multi-material design challenges. The problem setup is the same as depicted in \Cref{subfig:cantilever}, with the design domain's length-height ratio as $L:H=2:1$. For a more systematic assessment, we examined two distinct sets of five candidate materials, with the respective Young's moduli $E$ and normalized densities for the five candidate materials detailed in \Cref{tab:five-material-properties}.  
\begin{table}[ht]
    \centering
    \caption{Five-material cantilever: Properties of the five candidate materials.}
    \begin{tabular}{ccccccc}
        \toprule
        & MAT & 1 & 2 & 3 & 4 & 5 \\
        \midrule
        \multirow{3}{*}{Material Set 1} & $E_m$ & 0.4 & 0.7 & 0.85 & 0.9 & 1.0 \\
        & $\bar{M}_m$ & 0.3 & 0.5 & 0.65 & 0.8 & 1.0 \\
        & $E_m / \bar{M}_m$ & 1.33 & 1.4 & 1.31 & 1.13 & 1.0 \\
        \midrule
        \multirow{3}{*}{Material Set 2} & $E$ & 0.43 & 0.7 & 0.85 & 0.94 & 1.0 \\
        & $\bar{M}_m$ & 0.3 & 0.5 & 0.65 & 0.8 & 1.0 \\
        & $E_m / \bar{M}_m$ & 1.43 & 1.4 & 1.31 & 1.175 & 1.0 \\
        \bottomrule
    \end{tabular}
    \label{tab:five-material-properties}
\end{table}

\begin{figure*}[ht]
    \centering
    \begin{subfigure}[t]{.3\textwidth}
        \centering
        \includegraphics[width=0.9\textwidth]{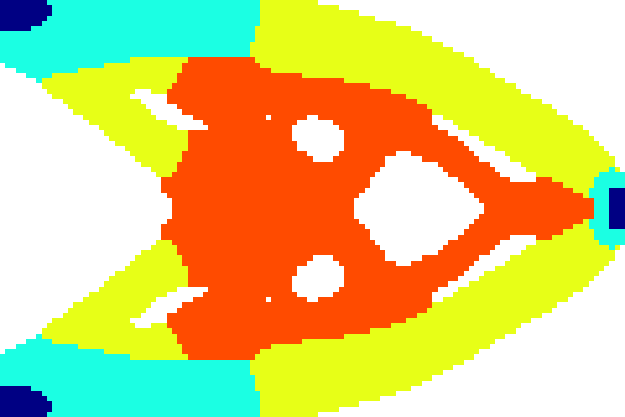}
        \caption{Discretization resolution: $120 \times 80$}
    \end{subfigure}
    \quad
    \begin{subfigure}[t]{.3\textwidth}
        \centering
        \includegraphics[width=0.9\textwidth]{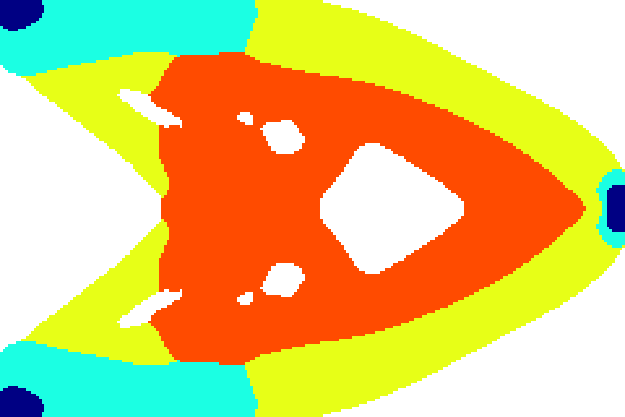}
        \caption{Discretization resolution: $240 \times 160$}
    \end{subfigure}
    \quad
    \begin{subfigure}[t]{.3\textwidth}
        \centering
        \includegraphics[width=0.9\textwidth]{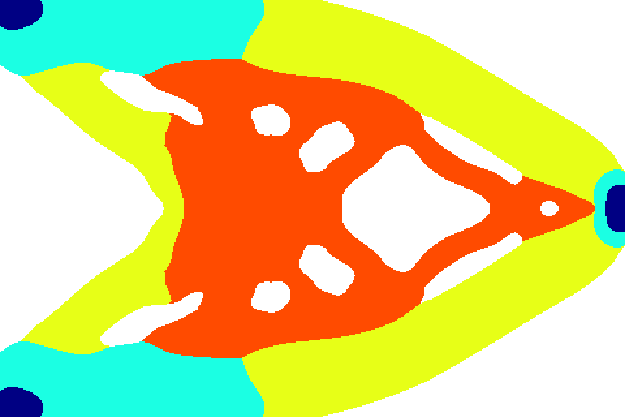}
        \caption{Discretization resolution: $480 \times 320$}
    \end{subfigure}
    \caption{Five-material cantilever: The resultant material configurations at various resolutions for Material Set 1. The color code for denoting different candidate materials is specified as: {\textcolor{mat1}{\rule{2ex}{2ex}}} denotes MAT 1;  {\textcolor{mat2}{\rule{2ex}{2ex}}} denotes MAT 2; {\textcolor{mat3}{\rule{2ex}{2ex}}} denotes MAT 3; {\textcolor{mat4}{\rule{2ex}{2ex}}} denotes MAT 4; and {\textcolor{mat5}{\rule{2ex}{2ex}}} denotes MAT 5.}
    \label{fig:five-material-minimum-compliance}
\end{figure*}

\begin{figure*}[ht]
    \centering
    \begin{subfigure}[t]{.3\textwidth}
        \centering
        \includegraphics[width=0.9\textwidth]{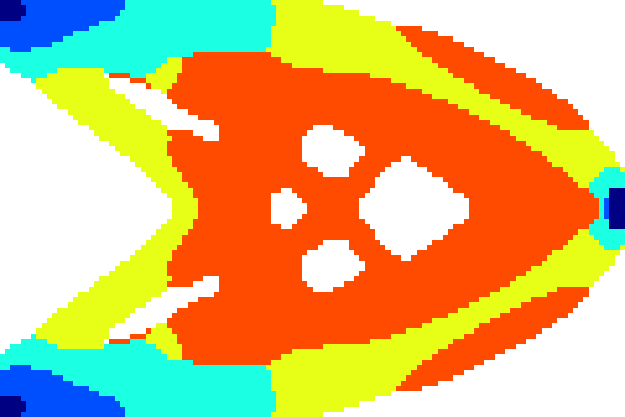}
        \caption{Discretization resolution: $120 \times 80$}
    \end{subfigure}
    \quad
    \begin{subfigure}[t]{.3\textwidth}
        \centering
        \includegraphics[width=0.9\textwidth]{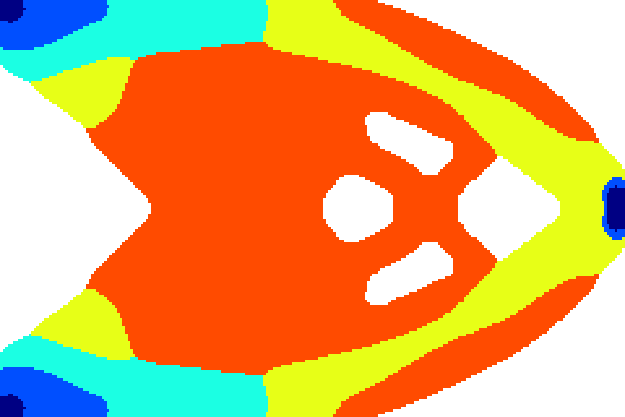}
        \caption{Discretization resolution: $240 \times 160$}
    \end{subfigure}
    \quad
    \begin{subfigure}[t]{.3\textwidth}
        \centering
        \includegraphics[width=0.9\textwidth]{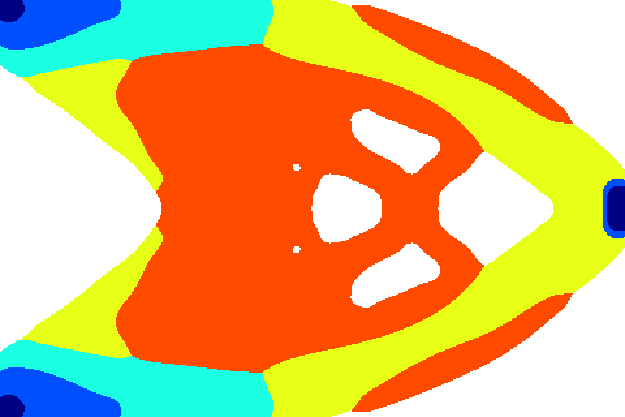}
        \caption{Discretization resolution: $480 \times 320$}
    \end{subfigure}
    \caption{Five-material cantilever: The resultant material configurations at various resolutions for Material Set 2. The color code for denoting different candidate materials is specified as: {\textcolor{mat1}{\rule{2ex}{2ex}}} denotes MAT 1;  {\textcolor{mat2}{\rule{2ex}{2ex}}} denotes MAT 2; {\textcolor{mat3}{\rule{2ex}{2ex}}} denotes MAT 3; {\textcolor{mat4}{\rule{2ex}{2ex}}} denotes MAT 4; and {\textcolor{mat5}{\rule{2ex}{2ex}}} denotes MAT 5.}
    \label{fig:five-material-minimum-compliance-2}
\end{figure*}

The parameter relaxation scheme was implemented with 8 distinct stages, as delineated in \Cref{tab:five-material-relaxation}, where the maximum permissible total mass fraction $\bar{M}_{\max}$ is decreased progressively from 0.5 to the designated target of 0.3 within the first $N_P = 7$ stages, following Eq. \eqref{eq:relaxation_interpolation}. In line with the previous tests, we maintained the same convergence tolerance of $\varepsilon = 5 \times 10^{-3}$.
\begin{table}[ht]
    \centering
    \caption{Five-material cantilever: Parameter relaxation scheme employed.}
    \begin{tabular}{ccccccccc}
        \toprule
        Stage & 1 & 2 & 3 & 4 & 5 & 6 & 7 & 8 \\
        \midrule
        $\bar{M}_{\max}$ & 0.500 & 0.459 & 0.422 & 0.387 & 0.356 & 0.327 & 0.300 & 0.300 \\
        $E_0$ & $10^{-2}$ & $10^{-2}$ & $10^{-2}$ & $10^{-2}$ & $10^{-2}$ & $10^{-2}$ & $10^{-2}$ & $10^{-9}$ \\
        \bottomrule
    \end{tabular}
    \label{tab:five-material-relaxation}
\end{table}
The results reported in the literature \cite{liu2024multi, huang2021new} are used to validate and compare with our results, as presented in \Cref{tab:five-material-minimum-compliance}. While the referenced papers \cite{liu2024multi, huang2021new} reported outcomes solely at a discretization resolution of $120 \times 80$, our numerical experiments extended to consider two higher resolutions as well.
\begin{table*}[ht]
    \centering
    \caption{Five-material cantilever: Main results, for both material sets considered and compared with the reference results reported in literature.}
    \label{tab:five-material-minimum-compliance}
    \resizebox{\textwidth}{!}{
    \begin{tabular}{ccccccccc}
        \toprule
        & \multirow{2}{*}{\textbf{Resolution}} & \multirow{2}{*}{$r$} & \multicolumn{2}{c}{\textbf{Our Method}} & \multicolumn{2}{c}{\textbf{Reference \cite{huang2021new}}} & \multicolumn{2}{c}{\textbf{Reference \cite{liu2024multi}}} \\
        \cmidrule(lr){4-5} \cmidrule(lr){6-7} \cmidrule(lr){8-9}
        & & & $N_{\text{FEM}}$ & $f$ & $N_{\text{FEM}}$ & $f$ & $N_{\text{FEM}}$ & $f$ \\
        \midrule
        \multirow{3}{*}{\textbf{Material Set 1}} & $120 \times 80$ & 3 & 34 & 37.500 & 250 & 37.881 & $>200$ & 36.568 \\
        & $240 \times 160$ & 6 & 36 & 37.724 & - & - & - & - \\
        & $480 \times 320$ & 12 & 36 & 38.458 & - & - & - & - \\
        \midrule
        \multirow{3}{*}{\textbf{Material Set 2}} & $120 \times 80$ & 3 & 33 & 36.631 & 250 & 36.722 & $>200$ & 35.755 \\
        & $240 \times 160$ & 6 & 35 & 36.614 & - & - & - & - \\
        & $480 \times 320$ & 12 & 34 & 37.325 & - & - & - & - \\
        \bottomrule
    \end{tabular}
    }
\end{table*}
The optimal material configurations obtained at various resolutions for two distinct material sets are illustrated in Figures. \ref{fig:five-material-minimum-compliance} and \ref{fig:five-material-minimum-compliance-2}.

As seen in \Cref{tab:five-material-minimum-compliance}, for both material sets our method yield comparable resultant objective function values with those reported in the state-of-the-art literature \cite{liu2024multi,huang2021new}. For Material Set 1, the obtained optimal material layout does not select all candidate materials but instead selects four out of five materials (i.e., MAT 1, 2, 3, and 5), as shown in \Cref{fig:five-material-minimum-compliance}. This is consistent with the finding reported in \cite{liu2024multi}, as the material preferences, implicitly embedded in the sensitivity specified in Eq. \eqref{eq:sensitivity_filtering_multimaterial}, align with those implied in \cite{liu2024multi}. For Material Set 2, since MAT 4 exhibits an increased Young's modulus ($E_4$) along with enhanced specific stiffness ($E_4/ \bar{M}_4$), all five candidate materials are retained in the final material layout, as depicted in \Cref{fig:five-material-minimum-compliance-2}, consistent with the findings reported in \cite{liu2024multi,huang2021new}.

By comparing the number of FEM analyses $N_{\text{FEM}}$, our method demonstrates its superior performance, achieving a reduction of approximately one order of magnitude in $N_{\text{FEM}}$. The iteration counts or $N_{\text{FEM}}$ were not provided in reference \cite{liu2024multi}. However, based on their reported number for the two-material case and the fact that the same methodology was employed for both the two-material and five-material cases, we estimate it to be more than 200. Furthermore, for both sets of candidate materials, the variance in $N_{\text{FEM}}$ across different discretization resolutions remains marginal, with a discrepancy of fewer than 3 iterations.  Significantly, increasing the discretization resolution does not correlate with an increased iteration count, suggesting that our framework maintains efficiency and scalability across various resolutions, demonstrating promise for tackling large-scale design problems. As the filter radius scales in accordance with the discretization resolution, the consistency of the resultant material layouts can be maintained by our proposed framework, as shown in Figures. \ref{fig:five-material-minimum-compliance} and \ref{fig:five-material-minimum-compliance-2}. That is, the primary structural features and the distribution of different candidate materials remain consistent.  Slight variations are mostly confined to the interfaces between distinct materials and some fine details, such as the configuration of minor structural elements like holes. This consistency across varying discretization resolutions, evidenced by both material sets, indicates the robustness of the proposed TO framework. Such robustness along with its computational efficiency and scalability affirm the superior capacity of our proposed TO framework to address multi-material design challenges.

\begin{figure*}[ht]
    \centering
    \begin{subfigure}[t]{.45\textwidth}
        \centering
        \includegraphics[width=0.75\textwidth]{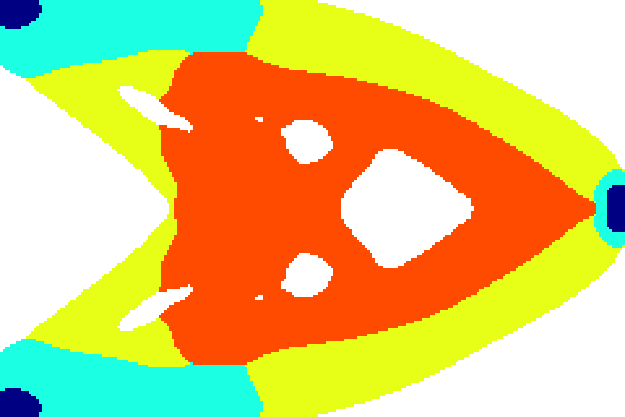}
        \caption{Material Set 1 (Disordered). $N_{\text{FEM}} = 34$, $f = 37.610$.}
    \end{subfigure}
    \quad \quad \quad
    \begin{subfigure}[t]{.45\textwidth}
        \centering
        \includegraphics[width=0.75\textwidth]{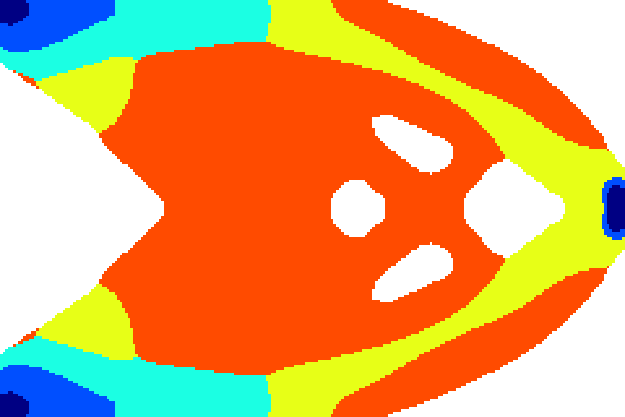}
        \caption{Material Set 2 (Disordered). $N_{\text{FEM}} = 35$, $f = 36.563$.}
    \end{subfigure}
    \caption{Five-material cantilever: The resultant material configurations  by solving the problems with the candidate materials disordered, at the discretization resolution of $240 \times 160$. The color code is the same as in \Cref{fig:five-material-minimum-compliance} and \Cref{fig:five-material-minimum-compliance-2}.}
    \label{fig:five-material-disordered}
\end{figure*}

Finally, we further validate that our methodology does not rely on any specific encoding of the candidate materials, i.e., it has no preference on the specific order of candidate materials, as the selection is solely based on the calculated sensitivity. To this end, the design variables $\boldsymbol{\rho}_1$-$\boldsymbol{\rho}_5$ are randomly ordered for representing MAT 1, MAT 3, MAT 5, MAT 4, MAT 2, respectively. With this new assignment of design variables, we re-solved the problem at the discretization resolution of $240 \times 160$. The results and optimal topologies are summarized in \Cref{fig:five-material-disordered}, where the materials are rendered with the original monolithic increasing order for ease of comparison. As compared with \Cref{fig:five-material-minimum-compliance}(b) and \Cref{fig:five-material-minimum-compliance-2}(b), the topologies are almost identical; the iteration counts and resultant objective function values are also similar to those reported in \Cref{tab:five-material-minimum-compliance}. These findings confirm that our proposed method does not need any specific sorting of candidate materials for multi-material TO designs.

\subsection{Multi-material Compliant Mechanism}

Finally, we investigated a compliant mechanism design incorporating five candidate materials, so as to assess the performance of our proposed framework for solving a non-convex TO problem with the complexity of involving multiple materials and mass constraints. Same as in the prior section, two distinct material sets, each with five candidate materials, were examined, with the respective Young's moduli $E$ and normalized densities for all candidate materials in each set outlined in \Cref{tab:five-material-properties}. The problem setup is the same as depicted in \Cref{fig:compliant_mechanism}. In the solution process, the parameter relaxation followed \textbf{Scheme 1} in \Cref{tab:two-material-relaxation}, i.e., a two-stage relaxation with $E_0 = 10^{-2}$ in the first stage and $E_0 = 10^{-9}$ in the second stage while keeping $\bar{M}_{\max} = 0.3$ in both stages. Being consistent with the other sections, the convergence tolerance was first set as $\varepsilon = 5 \times 10^{-3}$.  

The results are summarized in \Cref{tab:comparison_compliant_mechanism} and compared with those of the single-material case. The numbers of FEM analyses ($N_{\text{FEM}}$) required to reach convergence are consistently low (typically in twenties) across different discretization resolutions, same as the single-material compliant mechanism case and the five-material minimum compliance case, indicating the proposed framework's efficiency and scalability for dealing with multi-material non-convex designs. Due to the lack of literature data available for this test case, the results are not compared with any references. We note that the resultant displacement objective in this compliant mechanism design is improved by about 7\% across various discretization resolutions, by involving multiple candidate materials. While this difference may not be notably significant for the two material sets considered, the benefits of expanding material diversity in TO design can be magnified when considering factors such as cost, mass, and other material properties. 
\begingroup
\setlength{\tabcolsep}{10pt} 
\begin{table}[ht]
    \centering
    \caption{Five-material compliant mechanism: Main results and comparison with the single-material case.}
    \begin{tabular}{ccccc|cc}
        \toprule
        & \multirow{2}{*}{\textbf{Resolution}} & \multirow{2}{*}{$r$} & \multicolumn{2}{c}{\textbf{Five-material}} & \multicolumn{2}{c}{\textbf{Single-material}} \\
        \cmidrule(lr){4-5} \cmidrule(lr){6-7}
        & & & $N_{\text{FEM}}$ & $f$ & $N_{\text{FEM}}$ & $f$ \\
        \cmidrule(lr){1-7}
        \multirow{4}{*}{\textbf{Material Set 1}} & $200 \times 100$ & 2 & 22 & -0.9318 & 24 & -0.9226 \\
        & $400 \times 200$ & 4 & 22 & -0.8895 & 31 & -0.8602 \\
        & $600 \times 300$ & 6 & 23 & -0.8659 & 32 & -0.8413 \\
        & $800 \times 400$ & 8 & 31 & -0.8478 & 28 & -0.8167 \\
        \midrule
        \multirow{4}{*}{\textbf{Material Set 2}} & $200 \times 100$ & 2 & 22 & -0.9421 & 24 & -0.9226 \\
        & $400 \times 200$ & 4 & 25 & -0.8951 & 31 & -0.8602 \\
        & $600 \times 300$ & 6 & 24 & -0.8718 & 32 & -0.8413 \\
        & $800 \times 400$ & 8 & 29 & -0.8644 & 28 & -0.8167 \\
        \bottomrule
    \end{tabular}
    \label{tab:five-material-compliant_mechanism}
\end{table}
\endgroup

\begin{figure*}[ht]
    \centering
    \begin{subfigure}[t]{.3\textwidth}
        \centering
        \includegraphics[width=\textwidth]{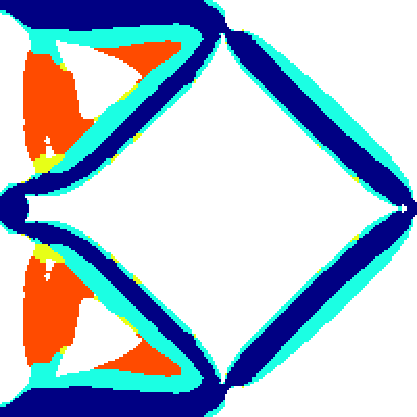}
        \caption{Discretization resolution $200 \times 100$.}
    \end{subfigure}
    \quad
    \begin{subfigure}[t]{.3\textwidth}
        \centering
        \includegraphics[width=\textwidth]{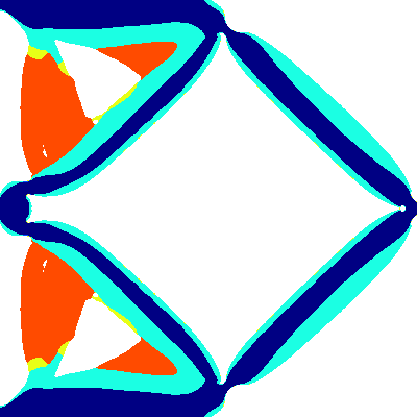}
        \caption{Discretization resolution $400 \times 200$.}
    \end{subfigure}
    \quad
    \begin{subfigure}[t]{.3\textwidth}
        \centering
        \includegraphics[width=\textwidth]{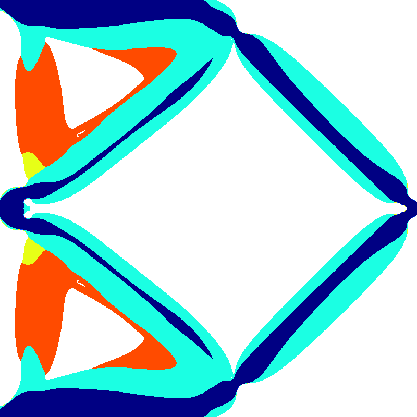}
        \caption{Discretization resolution $800 \times 400$.}
    \end{subfigure}
    \caption{Five-material compliant mechanism: The resulting material configurations at various resolutions for Material Set 1. The color code for denoting different candidate materials is specified as:  {\textcolor{mat1}{\rule{2ex}{2ex}}} denotes MAT 1;  {\textcolor{mat2}{\rule{2ex}{2ex}}} denotes MAT 2; {\textcolor{mat3}{\rule{2ex}{2ex}}} denotes MAT 3; {\textcolor{mat4}{\rule{2ex}{2ex}}} denotes MAT 4; and {\textcolor{mat5}{\rule{2ex}{2ex}}} denotes MAT 5.}
    \label{fig:five-material-compliant-mechanism}
\end{figure*}

\begin{figure*}[ht]
    \centering
    \begin{subfigure}[t]{.3\textwidth}
        \centering
        \includegraphics[width=\textwidth]{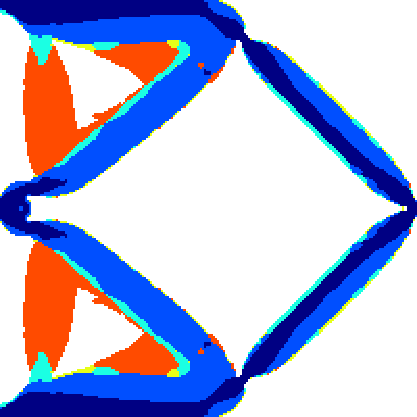}
        \caption{Discretization resolution $200 \times 100$.}
    \end{subfigure}
    \quad
    \begin{subfigure}[t]{.3\textwidth}
        \centering
        \includegraphics[width=\textwidth]{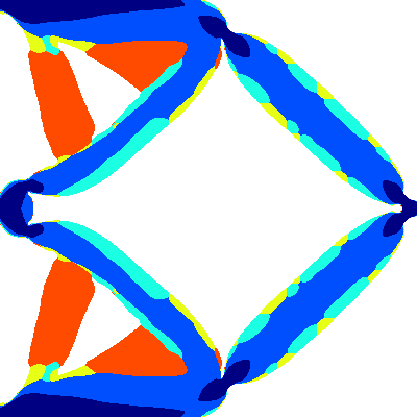}
        \caption{Discretization resolution $400 \times 200$.}
    \end{subfigure}
    \quad
    \begin{subfigure}[t]{.3\textwidth}
        \centering
        \includegraphics[width=\textwidth]{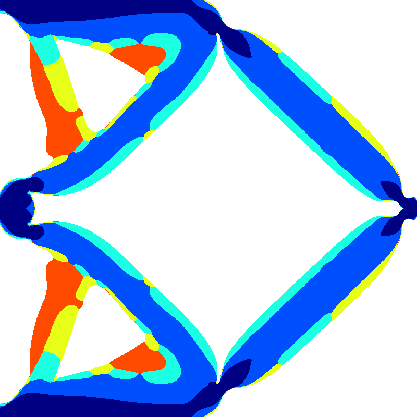}
        \caption{Discretization resolution $800 \times 400$.}
    \end{subfigure}
    \caption{Five-material compliant mechanism: The resulting material configurations at various resolutions for Material Set 2. The color code for denoting different candidate materials is specified as:  {\textcolor{mat1}{\rule{2ex}{2ex}}} denotes MAT 1;  {\textcolor{mat2}{\rule{2ex}{2ex}}} denotes MAT 2; {\textcolor{mat3}{\rule{2ex}{2ex}}} denotes MAT 3; {\textcolor{mat4}{\rule{2ex}{2ex}}} denotes MAT 4; and {\textcolor{mat5}{\rule{2ex}{2ex}}} denotes MAT 5.}
    \label{fig:five-material-compliant-mechanism-2}
\end{figure*}

The final material layouts are depicted in Figures \ref{fig:five-material-compliant-mechanism} and \ref{fig:five-material-compliant-mechanism-2}. Similarly to the five-material minimum compliance design, for Material Set 1, the design objective chose MAT 1, 2, 3 and 5 in the optimal layout; for Material Set 2, all five candidate materials were retained in the resultant material layout. For each material set, while lower resolutions yield rougher interfaces between different materials, the core structural attributes are maintained across various discretization resolutions.

\begin{figure*}[ht]
    \centering
    \begin{subfigure}[t]{.3\textwidth}
        \centering
        \includegraphics[width=\textwidth]{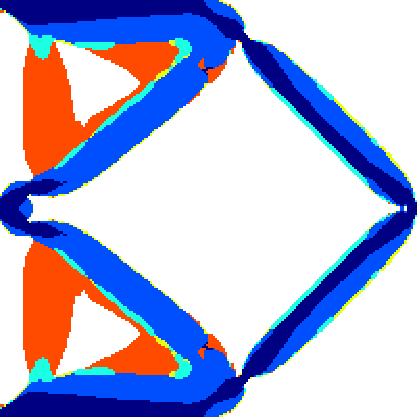}
        \caption{Discretization resolution $200 \times 100$. $\theta_1 = 0.7,~ \theta_2=1.3$.}
    \end{subfigure}
    \quad
    \begin{subfigure}[t]{.3\textwidth}
        \centering
        \includegraphics[width=\textwidth]{Figures/topo_cm_matset2_200x100_color.png}
        \caption{Discretization resolution $200 \times 100$. $\theta_1 = 0.7$, $\theta_2=1.5$.}
    \end{subfigure}
    \quad
    \begin{subfigure}[t]{.3\textwidth}
        \centering
        \includegraphics[width=\textwidth]{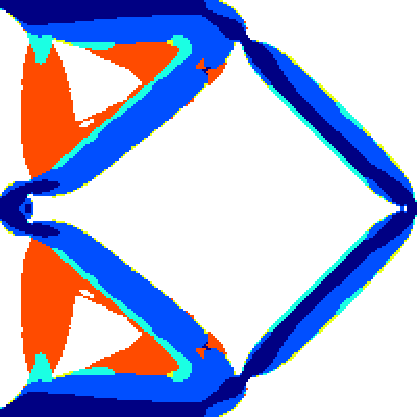}
        \caption{Discretization resolution $200 \times 100$. $\theta_1 = 0.7,~ \theta_2=1.7$.}
    \end{subfigure}
    \\
    \begin{subfigure}[t]{.3\textwidth}
        \centering
        \includegraphics[width=\textwidth]{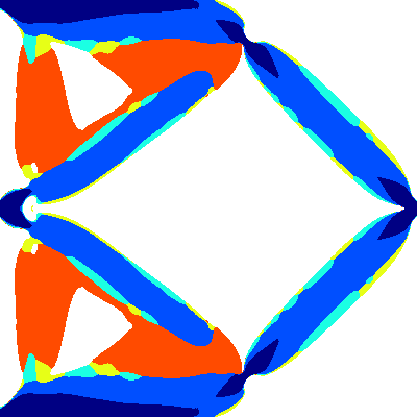}
        \caption{Discretization resolution $400 \times 200$. $\theta_1 = 0.6, ~\theta_2=1.5$.}
    \end{subfigure}
    \quad
    \begin{subfigure}[t]{.3\textwidth}
        \centering
        \includegraphics[width=\textwidth]{Figures/topo_cm_matset2_400x200_color.png}
        \caption{Discretization resolution $400 \times 200$. $\theta_1 = 0.7,~ \theta_2=1.5$.}
    \end{subfigure}
    \quad
    \begin{subfigure}[t]{.3\textwidth}
        \centering
        \includegraphics[width=\textwidth]{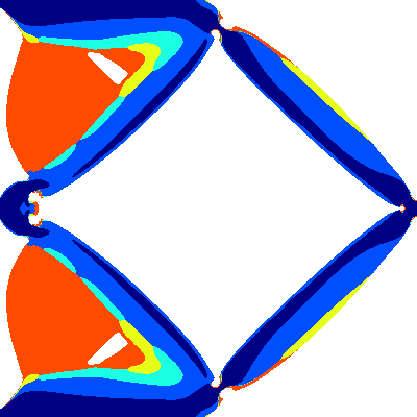}
        \caption{Discretization resolution $400 \times 200$. $\theta_1 = 0.8,~ \theta_2=1.5$.}
    \end{subfigure}
    \caption{Five-material compliant mechanism: The resulting material configurations for Material Set 2 with different adjustment factors ($\theta_1$ and $\theta_2$) for the trust-region radius. The color code for denoting different candidate materials is specified as:  {\textcolor{mat1}{\rule{2ex}{2ex}}} denotes MAT 1;  {\textcolor{mat2}{\rule{2ex}{2ex}}} denotes MAT 2; {\textcolor{mat3}{\rule{2ex}{2ex}}} denotes MAT 3; {\textcolor{mat4}{\rule{2ex}{2ex}}} denotes MAT 4; and {\textcolor{mat5}{\rule{2ex}{2ex}}} denotes MAT 5.}
    \label{fig:five-material-compliant-mechanism-adjustment-dependecy}
\end{figure*}

In the next test, we studied the effect of the two adjustment factors ($\theta_1$ and $\theta_2$) in Eq. \eqref{eq:trust_region_factor} for the trust-region radius. The values of $\theta_1$ and $\theta_2$ are varied with different combinations. As indicated in \Cref{tab:five-material-compliant_mechanism-adjustment-factor}, the optimization solutions show little dependency on these two factors: $N_{\text{FEM}}$ varies within $\pm 3$; the optimal objective function value $f$ varies within 6\%. As for the resulting topology, the two factors' effects are different,
as depicted in \Cref{fig:five-material-compliant-mechanism-adjustment-dependecy}. $\theta_2$, which controls the degree of enlarging the trust-region radius, shows minimum influence. However, $\theta_1$, the shrinking factor, has a more pronounced effect on the resultant topology. This is because shrinking the trust-region radius tends to increase the importance of the corresponding trust-region constraint in the multi-cut formulation \eqref{eq:multicut_trust_region}, and hence notably perturbs the resultant topology.

\textcolor{blue}{
\begin{table}[ht]
    \centering
    \caption{Five-material compliant mechanism: Results with different adjustment factors for the trust-region radius (i.e., $\theta_1$ and $\theta_2$ in Eq. \eqref{eq:trust_region_factor}).}
    \begin{tabular}{ccc|cc}
        \toprule
        \multirow{2}{*}{\textbf{Factor}} & \multicolumn{2}{c}{$200 \times 100$} & \multicolumn{2}{c}{$400 \times 200$} \\
        \cmidrule(lr){2-3} \cmidrule(lr){4-5}
        & $N_{\text{FEM}}$ & $f$ & $N_{\text{FEM}}$ & $f$ \\
        \cmidrule(lr){1-5}
        $\theta_1 = 0.7, \theta_2 = 1.5$ & 22 & -0.9421 & 25 & -0.8951 \\
        $\theta_1 = 0.7, \theta_2 = 1.7$ & 22 & -0.9333 & 27 & -0.9238 \\
        $\theta_1 = 0.7, \theta_2 = 1.3$ & 22 & -0.9351 & 25 & -0.9032 \\
        $\theta_1 = 0.6, \theta_2 = 1.5$ & 25 & -0.9608 & 25 & -0.9107 \\
        $\theta_1 = 0.8, \theta_2 = 1.5$ & 24 & -0.9466 & 28 & -0.8729 \\
        \bottomrule
    \end{tabular}
    \label{tab:five-material-compliant_mechanism-adjustment-factor}
\end{table}
}

\begin{figure*}[ht]
    \centering
    \begin{subfigure}[t]{.3\textwidth}
        \centering
        \includegraphics[width=\textwidth]{Figures/topo_cm_matset2_800x400_color.png}
        \caption{Convergence tolerance $5 \times 10^{-3}$.}
    \end{subfigure}
    \quad
    \begin{subfigure}[t]{.3\textwidth}
        \centering
        \includegraphics[width=\textwidth]{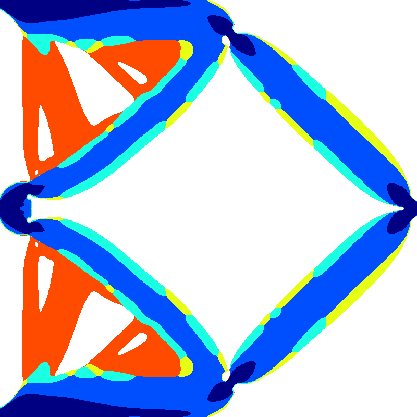}
        \caption{Convergence tolerance $10^{-3}$.}
    \end{subfigure}
    \quad
    \begin{subfigure}[t]{.3\textwidth}
        \centering
        \includegraphics[width=\textwidth]{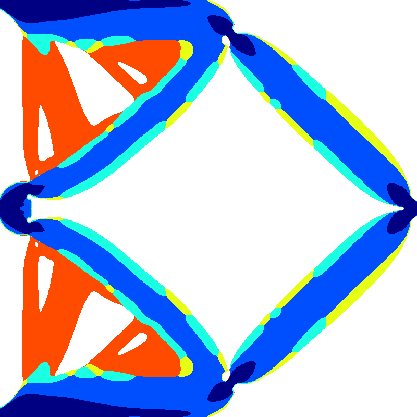}
        \caption{Convergence tolerance $5 \times 10^{-4}$.}
    \end{subfigure}
    \caption{Five-material compliant mechanism: The resulting material configurations obtained with different convergence tolerance $\varepsilon$, all at the discretization resolution of $800 \times 400$ and for Material Set 2. The color code for denoting different candidate materials is specified as:  {\textcolor{mat1}{\rule{2ex}{2ex}}} denotes MAT 1;  {\textcolor{mat2}{\rule{2ex}{2ex}}} denotes MAT 2; {\textcolor{mat3}{\rule{2ex}{2ex}}} denotes MAT 3; {\textcolor{mat4}{\rule{2ex}{2ex}}} denotes MAT 4; and {\textcolor{mat5}{\rule{2ex}{2ex}}} denotes MAT 5.}
    \label{fig:five-material-compliant-tolerance}
\end{figure*}

In the last test, we systematically studied the impact of the convergence tolerance $\varepsilon$ set in \Cref{algo:multi-cuts} on the final solution. Specifically, we examined three different tolerances, as outlined in \Cref{tab:five-material-compliant_mechanism}. By comparing $N_{\text{FEM}}$ and the optimal objective function value $f$, we find that decreasing the tolerance leads to more optimization iterations to reach a better solution (or a lower objective function value). However, once the best solution has been found, further decreasing the tolerance does not change the solution, as indicted in both \Cref{tab:five-material-compliant_mechanism-tolerance} and \Cref{fig:five-material-compliant-tolerance}. And it is worth noting that our method did not result in redundant iterations or fluctuating objective function values when further decreasing the tolerance, which contrasts with the behavior observed in other methods, e.g., \cite{huang2010evolutionary, liang2019topology}.
\begingroup
\setlength{\tabcolsep}{10pt} 
\begin{table}[ht]
    \centering
    \caption{Five-material compliant mechanism: Results with different convergence tolerance $\varepsilon$.}
    \begin{tabular}{cccc|cc|cc}
        \toprule
        \multirow{2}{*}{\textbf{Resolution}} & \multirow{2}{*}{$r$} & \multicolumn{2}{c}{$\varepsilon = 5 \times 10^{-3}$} & \multicolumn{2}{c}{$\varepsilon = 10^{-3}$} & \multicolumn{2}{c}{$\varepsilon = 5 \times 10^{-4}$} \\
        \cmidrule(lr){3-4} \cmidrule(lr){5-6} \cmidrule(lr){7-8}
        & & $N_{\text{FEM}}$ & $f$ & $N_{\text{FEM}}$ & $f$ & $N_{\text{FEM}}$ & $f$ \\
        \cmidrule(lr){1-8}
        $200 \times 100$ & 2 & 22 & -0.9421 & 36 & -0.9738 & 32 & -0.9738 \\
        $400 \times 200$ & 4 & 25 & -0.8951 & 39 & -0.9348 & 39 & -0.9348 \\
        $600 \times 300$ & 6 & 24 & -0.8718 & 24 & -0.8718 & 24 & -0.8718 \\
        $800 \times 400$ & 8 & 29 & -0.8644 & 35 & -0.8770 & 35 & -0.8770 \\
        \bottomrule
    \end{tabular}
    \label{tab:five-material-compliant_mechanism-tolerance}
\end{table}
\endgroup
\section{Conclusion}
\label{sec:conclusion}
We have presented a new TO framework that can efficiently solve both convex and non-convex TO problems involving single or multiple materials. It directly handles binary design variables, eliminating the need for additional efforts to convert binary variables into continuous variables and vice versa (e.g., through interpolation and projection as in SIMP). When dealing with multi-material TO designs, it does not need to enumerate all possible combinations of candidate materials and solve multiple independent TO problems like SIMP \cite{zuo2017multi}, nor does it need to sort candidate materials in a specific order and interpolate design variables between different materials like other methods \cite{zuo2017multi, huang2021new, liu2024multi}. Instead, it only needs to expand the dimension of design variables and incorporate additional inequality constraints arising from the mass constraint. 

Through numerical tests and comparisons with other methods, we have demonstrated that our framework can significantly reduce the number of optimization iterations (and thereby the number of FEM analyses required) by about one order of magnitude, while maintaining comparable solution quality in terms of the optimal value of the objective function and the resulting material layout. Despite increasing discretization resolutions and the inclusion of multiple materials—which significantly increase the number of design variables and introduce more inequality constraints—we observe consistency in solution quality and the number of optimization iterations required to reach the optimal solution. We have also found that the minimum Young's modulus has the greatest impact on the conditioning of the optimization problem and is therefore considered as the primary parameter for relaxation. Relaxing it, along with a secondary parameter in multiple stages, could potentially offer greater advantages in achieving lower objective function values. The extent of improvement can vary depending on the specific problem, particularly its conditioning.

Looking ahead, we envision that the new TO framework proposed in this work will be particularly beneficial for large-scale TO applications, which entail numerous design variables and constraints and demand substantial computational resources for FEM analyses (or PDE solving). Additionally, the inclusion of nonlinear governing equations \cite{picelli2022topology} or nonlinear material models \cite{wang2014interpolation} represents a natural extension, broadening the framework's utility and relevance for more complex TO designs. Furthermore, incorporating a broader class of constraints into this framework, such as stress constraints \cite{kiyono2023stress} and manufacturability constraints \cite{sigmund2009manufacturing}, can be another interesting future direction. Finally, in this work, the single-cut problem in Eq. \eqref{eq:single-cut} and the multi-cut problems in Eq. \eqref{eq:branched_multi-cut} were all solved using an off-the-shelf mixed integer linear programming (MILP) solver \cite{gurobi}. It could potentially be replaced by a more advanced MILP solution method with improved efficiency, which would be particularly beneficial when dealing with multi-material TO. In such cases, the MILP problems with substantially expanded binary variables require a faster optimizer. One possible direction is to explore the geometric multi-grid characteristics inherent in material layouts and leverage these properties to accelerate optimizers by utilizing the transferability of Graph Neural Networks (GNNs) across optimization problems of varying sizes \cite{cappart2023combinatorial}.

\bibliographystyle{elsarticle-num}
\bibliography{ref}
\end{document}